\documentclass[preprint]{elsarticle}

\pdfoutput=1

\usepackage[left=3cm,right=3cm]{geometry}
\usepackage{amsmath}
\usepackage{amsthm}
\usepackage{amssymb}
\usepackage{graphicx}
\usepackage{algorithmic}
\usepackage{algorithm}
\usepackage{verbatim}
\usepackage{color}
\usepackage{setspace}

\def\q{\mathbf{q}}
\def\Q{\mathbf{Q}}
\def\M{\mathbf{M}}
\def\etal{et al.}
\def\ie{i.e.}

\newcommand{\unit}[1]{\ensuremath{\, \mathrm{#1}}}

\title{Asynchronous Variational Contact Mechanics}
\author[CU]{E. Vouga\corref{cor}}
\author[CU]{D. Harmon}
\author[DISNEY]{R. Tamstorf}
\author[CU]{E. Grinspun}

\cortext[cor]{Corresponding author}

\address[CU]{Computer Science Department, Columbia University, New York, New York 10027}
\address[DISNEY]{Walt Disney Animation Studios, Burbank, California 91521}

\begin{document}

\begin{abstract}
An asynchronous, variational method for simulating elastica in complex contact and impact scenarios
is developed. Asynchronous Variational Integrators~\cite{aLew2003} (AVIs) are extended to handle contact forces by associating different time steps to forces instead of to spatial elements. By discretizing a barrier potential by an infinite sum of nested quadratic potentials, these extended AVIs are used to resolve contact while
obeying momentum- and energy-conservation laws. A series of two- and three-dimensional
examples illustrate the robustness and good energy behavior of the method.
\end{abstract}


\maketitle

\section{Introduction}

\emph{Variational integrators} (VIs)~\cite{aSuris1990, cMacKay1992,
  aMarsden2001} are a general class of time integration methods for
Hamiltonian systems whose construction guarantees certain properties highly
desirable of numerical simulations. Instead of directly discretizing the smooth
equations of motion of a system, the variational approach instead
discretizes the system's action integral. By analogy to Hamilton's least
action principle, a discrete action can be formed, and \emph{discrete}
Euler-Lagrange equations derived by examining paths which extremize
it. From the Euler-Lagrange equations, discrete equations of motion
are readily recovered. As a consequence of this special construction, VIs are guaranteed to satisfy a discrete
formulation of Noether's Theorem~\cite{tWest2004}, and as a special
case conserve linear and angular momentum. VIs are automatically
\emph{symplectic}~\cite{bHairer2006}; while they do not necessarily
conserve energy, conservation of the symplectic form assures no-drift
conservation of energy over exponentially many time
steps~\cite{bHairer2006}.

Given the many advantages of VIs, it is natural to apply them to the handling of contact and impact, a long-studied and
challenging problem in physical simulation. Unfortunately, a na{\"\i}ve
application of a contact algorithm to a variational integrator is not
guaranteed to preserve the variational structure of the time
integration method, and in practice one observes that the good energy
behavior is lost. For this reason, a few recent works have explored
structure-preserving approaches for contact
mechanics~\cite{aBarth1999,aKane2000,aLaursen2002,aPandolfi2002}.
Common to all these approaches is a \emph{synchronous} treatment of
global time, in which the entire configuration is advanced from one
intant in time to the next. While synchronous integration is
attractive for its simplicity, it has the drawback that a
spatially-localized stiff mode---such as that associated to a
localized contact---can force the global configuration to
advance at fine time steps. 

Indeed, mechanical systems are almost never uniformly stiff. Different
potentials have different stable time step requirements, and even for
identical potentials this requirement depends on element size, since
finer elements can support higher-energy modes than coarser
elements. Any global time integration scheme cannot take advantage of
this variability, and instead must integrate the entire system at the
globally stiffest time step. Suppose the system can be partitioned
into elements such that each force acts entirely within one
element. Then \emph{asynchronous variational integrators} (AVIs)~
\cite{aLew2003} generalize VIs by allowing each element to have its
own, independent time step. Coarser elements can then be assigned a
slower ``clock,'' and finer elements a faster one. Asynchrony avoids
the undesirable situation in which a small number of very fine
elements degrade overall performance.  AVIs retain all of the
properties of variational integrators mentioned above, except for
discrete symplecticity. However, AVIs instead preserve an analogous
\emph{discrete multisymplectic} form, and it has been shown experimentally that
preservation of this form likely induces the same long-time good
energy behavior that characterize symplectic
integrators~\cite{aLew2003}. 

To our knowledge, this work is the first to consider an asynchronous,
variational treatment of contact shown to retain multisymplecity. We are also aware that Ryckman and Lew~\cite{aLew201X} 
are concurrently investigating extending the AVI framework to incorporate contact response.

The starting point for this approach is the selection of the penalty
method as a model for contact~\cite{bWriggers2002, bLaursen2002}. For each pair of
elements in the system, a potential is added
that is (piecewise) quadratic in the \emph{gap function} measuring the separation distance between the
two elements. This potential vanishes when elements are sufficiently
far apart, and increases with increasing interpenetration, so
that approaching elements feel a force that resists impact. This
approach suffers two limitations, however. Firstly, these contact
potentials are fundamentally nonlocal phenomena: for every pair of
elements that might come into contact during the course of the
simulation, a potential coupling the two must be added. As will be
shown, the fact that contact potentials cannot be expressed as the
integration over the material domain of an energy density depending only on a neighborhood of the domain will present a technical
obstruction to the original formulation of AVIs, but fortunately one
that can be overcome by a natural generalization.

Secondly, penalty forces have a well-studied performance-robustness
tradeoff~\cite{aBelytschko1991}: adding a half-quadratic potential requires
choosing an arbitrary stiffness parameter, and for any stiffness
chosen for the penalty potential, two approaching elements will interpenetrate some distance, and in the worst case \emph{tunnel} completely through
each other. Moreover, the stable time step of the penalty force
decreases as stiffness increases, so choosing a very stiff penalty
potential is untenable as a solution to excessive penetration or tunneling. In practice, users of the method must determine
an adequate penalty stiffness by iterated
tweaking of parameters, until the simulation completes without
collision artifacts. An appealing modification of the penalty approach replaces the quadratic potential with a nonlinear \emph{barrier potential}~\cite{bNocedal2000} that diverges as the configuration
approaches contact. Because the barrier diverges, its stiffness is
unbounded, necessitating a time-adaptive time stepping method. This work
presents a discrete analogue of the barrier potential---an infinite
sequence of \emph{discrete penalty layers}---that in effect enables
AVIs to serve as adaptive integrators.

This paper
\begin{itemize}
\item{extends the construction of AVIs so that a discretization into
    disjoint elements is no longer necessary, by associating a clock
    to each force instead of to each element;}

\item{demonstrates that this generalization does not destroy the
    desirable integration properties guaranteed by the variational
    paradigm, most importantly the conservation of the discrete multisymplectic
    form;}

\item{leverages this extension to equip the AVI framework with a
    contact model. The proposed barrier method uses a divergent
    sequence of quadratic potentials that guarantees non-penetration
    and retains the asynchrony or conservation properties of AVIs;}

\item{presents numerical evidence to support the claim that by
    retaining the symplectic structure of the smooth system,
    simulations of thin shells undergoing complex (self-)interactions
    have demonstrably good long-time energy and momentum behavior;}

\item{describes simple extensions to the contact model to allow for
    controlled, dissipative phenomena, such as a coefficient of
    restitution and kinetic friction. Although there is not yet theory
    explaining the energy behavior of dissipative simulations run
    under a variational integrator, empirical evidence is presented to show that the
    proposed method produces smooth, controlled, and qualitatively
    correct energy decay.}
\end{itemize}

This paper complements the publication~\cite{aHarmon2009}, which
provides a detailed description of the software implementation using
\emph{kinetic data structures}~\cite{aGuibas1998}. For completeness,
Section \ref{sKDS} briefly introduces these concepts.

\section{Related Work}
The simplest contact models for finite element simulation follow the
early analytical work of Hertz~\cite{aHertz1882} in assuming
frictionless contact of planar (or nearly planar) surfaces with small
strain. In this regime, several approaches have been explored to
arrive at a weak formulation of contact; for a high-level survey of
these approaches, see for example the overview by Belytschko
\etal~\cite{bBelytschko2006} or Wriggers~\cite{aWriggers1995}. The
first of these are the use of penalty forces, described for instance by
Oden~\cite{cOden1980} and Kikuchi and Oden~\cite{bKikuchi1988}. The
penalty approach results in a contact force proportional to an arbitrary
\emph{penalty stiffness} parameter and to the rate of interpenetration, or
in more general formulations to an arbitrary function of rate of
interpenetration and interpenetration depth; Belytschko and
Neal~\cite{aBelytschko1991} discuss the choosing of this parameter in
Section 8. Recent work by Belytschko \etal~\cite{aBelytschko2002} uses
moving least squares to construct an implicit smooth contact surface,
from which the interpenetration distance is evaluated. Peric and
Owen~\cite{aPeric1992} describe how to equip penalty forces with a
Coulomb friction model.

Seeking to exactly enforce non-penetration along the contact surface
leads to generalizations of the method of Lagrange multipliers. Hughes
\etal~\cite{aHughes1976} and Nour-Omid and
Wriggers~\cite{aNourOmid1986} provide an overview of this approach in
the context of contact response. Such contraint enforcement can be
viewed as a penalty force in the limit of infinite stiffness,
impossible to attain in practice since the system becomes
ill-conditioned. Taylor and Papadopoulos~\cite{aTaylor1993} considers
persistent contact by extending Newmark to treat jump conditions in
kinematic fields, thus reducing undesirable oscillatory
modes. However, the effects of these modifications on numerical
dissipation and long-time energy behavior is not considered.

The Augmented Lagragian method blends the penalty and Lagrange multiplier
approaches, and combines the advantages of both: unlike for pure
penalty forces, convergence to the exact interpenetration constraint
does not require taking the penalty stiffness to infinity, and the
Lagrange multiplier solve tends to be well-conditioned. Bertsekas
\cite{bBertsekas1984} gives a mathematical overview of the augmented
Lagrangian method, and Wriggers \etal~\cite{aWriggers1985} and Simo
and Laursen~\cite{aSimo1992} expand on its application to contact
problems in finite elements.

Non-smooth contact requires special consideration, since in the
non-smooth regime there is no straightforward way of defining a
contact normal or penetration distance. Simo \etal~\cite{aSimo1985}
discretize the contact surface into segments over which they assume
constant contact pressure; this formulation allows them to handle
non-node-to-node contact using a perturbed Lagrangian. Kane
\etal~\cite{aKane1999} apply non-smooth analysis to resolve contact
constraints between sharp objects. Pandolfi \etal~\cite{aPandolfi2002}
extend the work of Kane \etal\,by describing a variational model for
non-smooth contact with friction. Cirak and West~\cite{aCirak2005}
decompose contact resolution into an impenetrability-enforcement and
momentum-transfer step, thereby exactly enforcing non-interpenetration
while nearly conserving momentum and energy.

Several authors have explored a structure-preserving approach to
solving the contact problem. Barth \etal~\cite{aBarth1999} consider an
adaptive-step-size algorithm that preserves the time-reversible
symmetry of the RATTLE algorithm, and demonstrate an application to
an elastic rod interacting with a Lennard-Jones potential. Kane
\etal~\cite{aKane2000} show that the Newmark method, for all
parameters, is variational, and construct two two-step dissipative
integrators that yield good energy decay. Laursen and Love~\cite{aLaursen2002}, by taking into account
velocity discontinuities that occur at contact interfaces, develop 
a momentum- and energy-preserving method for simulating frictionless contact. This paper shares with these last
approaches the viewpoint that structured integration, with its associated
conservation guarantees, is an invaluable tool for accurately simulating dynamic 
systems with contact. 

Although several previous approaches are also adaptive, the
algorithm described in this paper is the first structured integrator for contact mechanics
that achieves time adaptivity using asynchrony. This novel approach guarantees the robustness
of the proposed integrator, without compromising the good properties of structured
integration.

\section{Variational Integrators}
This section presents a background on variational integration and symplectic structure~\cite{bHairer2006, aMarsden2001, tWest2004}.

Let $\gamma(t)$ be a piecewise-regular trajectory through 
configuration space $\Q$, and $\dot\gamma(t) = \frac{d}{dt}\gamma(t)$ 
be the configurational velocity at time $t$. For simplicity, 
assume that the kinetic energy of the system $T$ depends only on 
configurational velocity, and that the potential energy $V$ depends 
only on configurational position, so that the Lagrangian 
$L$ at time $t$ may be written as
\begin{align}
L(q, \dot q) = T(\dot q) - V(q). \label{eq-action}
\end{align}

Then given the configuration of the system $q_0$ at time $t_0$ and $q_f$ at $t_f$, Hamilton's principle~\cite{bLanczos1986} states that the trajectory of the system $\gamma(t)$ joining $\gamma(t_0) = q_0$ and $\gamma(t_f) = q_f$ is a stationary point of the action functional
\begin{align*}
S(\gamma) = \int_{t_0}^{t_f} L\left[ \gamma(t), \dot\gamma(t) \right] dt
\end{align*}
with respect to taking variations $\delta \gamma$ of $\gamma$ which leave $\gamma$ fixed at the endpoints $t_0,\, t_f$. In other words, $\gamma$ satisfies
\begin{align}
dS(\gamma) \cdot \delta \gamma = 0. \label{eq-Hamilton}
\end{align}

Integrating by parts, and using that $\delta \gamma$ vanishes at $t_0$ and $t_1$,
\begin{align*}
dS(\gamma) \cdot \delta \gamma = \int_{t_0}^{t_f} \left( \frac{\partial L}{\partial q}(\gamma,\dot\gamma) \cdot \delta \gamma + \frac{\partial L}{\partial \dot q} (\gamma, \dot \gamma) \cdot \delta \dot\gamma \right) dt = \int_{t_0}^{t_f} \left( -\frac{\partial V}{\partial q}(\gamma) - \frac{\partial^2 T}{\partial \dot q^2} (\dot \gamma)\ddot\gamma \right) \cdot \delta \gamma \, dt = 0.
\end{align*}

Since this equality must hold for all variations $\delta \gamma$ that fix $\gamma$'s endpoints, 
\begin{align}
\frac{\partial V}{\partial q}(\gamma) + \frac{d}{dt}\left(\frac{\partial T}{\partial \dot q}(\dot \gamma)\right) = 0 \label{eq-EL},
\end{align}
the \emph{Euler-Lagrange equation} of the system. This equation is a second-order ordinary differential equation, and so has a unique solution $\gamma$ given two initial values $\gamma(t_0)$ and $\dot \gamma(t_0)$.

\subsection{Symplecticity}

The flow $\Theta_s: \left[\gamma(t), \dot \gamma(t) \right] \mapsto \left[\gamma(t+s), \dot \gamma(t+s)\right]$ induced by (\ref{eq-EL}) has many structure-preserving properties; in particular it is momentum-preserving, energy-preserving, and symplectic~\cite{tLew2003}. To derive this last property, for the remainder of this section the space of trajectories is restricted to those that satisfy the Euler-Lagrange equations. For such trajectories, if the requirement that $\delta\gamma$ fix the endpoints of $\gamma$ is relaxed, then the boundary terms of the integration by parts are no longer 0 and
\begin{align}
dS(\gamma) \cdot \delta \gamma = \frac{\partial T}{\partial \dot q}\left[ \pi_{\dot q} (q, \dot q) \right]\cdot \delta \gamma\bigg|_{t_0}^{t_f},\label{eq-ELparam}
\end{align}
where $\pi_{\dot q}$ is projection onto the second factor.

Since initial conditions $(q, \dot q)$ are in bijection with 
trajectories satisfying the Euler-Lagrange equation, such trajectories 
$\gamma$ can be uniquely parametrized by initial conditions 
$\left[\gamma(t_0), \dot\gamma(t_0)\right]$. For the remainder of this 
section variations $\delta \gamma$ are also restricted to \emph{first variations}: those variations in whose direction
$\gamma$ continues to satisfy the Euler-Lagrange equations. These 
are also parametrized by variations of the initial conditions, 
$(\delta q, \delta \dot q)$. For conciseness of notation, the change of variables
$\nu(t) = (\gamma(t), \dot \gamma(t))$ and $\delta \nu(t) = 
\left[\delta \gamma(t), \delta \dot \gamma(t)\right]$ can be used; using this 
notation the above two facts can be rewritten as $\nu(t) = 
\Theta_{t-t_0}\nu(t_0)$ and $\delta\nu(t) = {\Theta_{t-t_0}}_* 
\delta\nu(t_0)$. The action (\ref{eq-action}), a functional on 
trajectories $\gamma$, can also be rewritten as a function $S_i$ of 
the intial conditions,
\begin{align*}
S_i(q, \dot q) = \int_0^{t_f-t_0} L\left[\Theta_t(q, \dot q) \right]dt,
\intertext{so that}
dS(\gamma) \cdot \delta \gamma = dS_i\left[\nu(t_0) \right] \cdot \delta \nu(t_0).
\intertext{Substituting all of these expressions into (\ref{eq-ELparam}),}
dS_i\left[\nu(t_0)\right] \cdot \delta \nu(t_0) &= \left(\frac{\partial T}{\partial \dot q} \circ \pi_{\dot q}\right) \left[\Theta_{t-t_0}\nu(t_0)\right] \cdot \delta \gamma(t)\bigg|_{t_0}^{t_f}\\
&= \left(\frac{\partial T}{\partial \dot q} \circ \pi_{\dot q} \right) \left[\Theta_{t-t_0}\nu(t_0)\right] dq \cdot \delta \nu(t)\bigg|_{t_0}^{t_f}\\
&= \left(\frac{\partial T}{\partial \dot q} \circ \pi_{\dot q} \right) \left[\Theta_{t-t_0}\nu(t_0)\right] dq \cdot {\Theta_{t-t_0}}_* \delta \nu(t_0)\bigg|_{t_0}^{t_f}\\
&= ({\Theta_{t_f-t_0}}^* \theta_L - \theta_L)_{\nu(t_0)} \cdot \delta \nu(t_0),
\intertext{where $\theta_L$ is the one-form $\left(\frac{\partial T}{\partial \dot q}\circ \pi_{\dot q}\right) dq$. Since $dS_i$ is exact,}
d^2S_i &= 0 = {\Theta_{t_f-t_0}}^*d\theta_L - d\theta_L,
\end{align*}
so since $t_0$ and $t_f$ are arbitrary, $\Theta_s^*d\theta_L = d\theta_L$ for arbitrary times $s$, and $\Theta$ preserves the so-called \emph{symplectic form} $d\theta_L$.

\subsection{Discretization \label{s-discrete}}

Discrete mechanics~\cite{aVeselov1988, aSuris1990, aMoser1991, aWendlandt1997, aMarsden2001, bHairer2006} describes a discretization of Hamilton's principle, yielding a numerical integrator that shares many of the structure-preserving properties of the continuous flow $\Theta_s$. Consider a discretization of the trajectory $\gamma: [t_0, t_f] \to \Q$ by a piecewise linear trajectory interpolating $n$ points $\q = \{q_0, q_1, \ldots q_{n-1}\}$, with $q_0 = \gamma(t_0)$ and $q_{n-1} = \gamma(t_f)$, where the discrete velocity $\dot q_{i+1/2}$ on the segment between $q_i$ and $q_{i+1}$ is
\begin{align*}
\dot q_{i+1/2} = \frac{q_{i+1}-q_{i}}{h}, \quad h = \frac{t_f-t_0}{n}.
\end{align*}

An analogue of (\ref{eq-EL}) in this discrete setting is needed. To that end, a discrete Lagrangian
\begin{align}
L_d(q_a, q_b) = T\left(\frac{q_b-q_a}{h}\right) - V(q_b)\label{eq-DL}
\intertext{can be formulated, as well as a discrete action}
S_d(\q) = \sum_{i=0}^{n-2} h L_d(q_i, q_{i+1}). \label{eq-daction}
\end{align}

Motivated by (\ref{eq-Hamilton}), a discrete Hamilton's principle can be imposed:
\begin{align*}
dS_d(\q) \cdot \delta \q = 0
\end{align*}
for all variations $\delta \q = \{\delta q_0, \delta q_1, \ldots, \delta q_{n-1}\}$ that fix $\q$ at its endpoints, \ie, with $\delta q_0 = \delta q_{n-1} = 0.$ For ease of notation, the kinetic and potential energy terms in (\ref{eq-DL}) can be written to depend on $(q_a, q_b)$, two points of phase space consecutive in time, instead of $(q, \dot q)$:
\begin{align*}
T_d(q_a, q_b) &= T\left(\frac{q_b-q_a}{h}\right) &
T'_d(q_a, q_b) &= \frac{\partial T}{\partial \dot q}\left(\frac{q_b-q_a}{h}\right)\\
V_d(q_a, q_b) &= V(q_b) &
V'_d(q_a, q_b) &= \frac{\partial V}{\partial q}(q_b).
\end{align*}
Then
\begin{align*}
dS_d(\q) \cdot \delta \q &= \sum_{i=0}^{n-2} h \left(D_1 L_d(q_i, q_{i+1}) \cdot \delta q_i + D_2 L_d(q_i, q_{i+1}) \cdot \delta q_{i+1}\right)\\
&= \sum_{i=0}^{n-2} h \left(-\frac{1}{h}T'_d(q_i,q_{i+1}) \cdot \delta q_i + \frac{1}{h}T'_d(q_i, q_{i+1}) \cdot \delta q_{i+1} - \frac{\partial V}{\partial q}(q_{i+1}) \cdot \delta q_{i+1} \right)\\
&= T'_d(q_{n-2},q_{n-1})\cdot \delta q_{n-1} - T'_d(q_0,q_1)\cdot \delta q_0 - h \frac{\partial V}{\partial q}(q_{n-1}) \cdot \delta q_{n-1}\\
&\quad + \sum_{i=1}^{n-2} \left( T'_d(q_{i-1},q_i) - T'_d(q_i,q_{i+1}) - h \frac{\partial V}{\partial q}(q_i)\right)\cdot \delta q_i \\
&= \sum_{i=1}^{n-2} \left( T'_d(q_{i-1}, q_{i}) - T'_d(q_i, q_{i+1}) - h \frac{\partial V}{\partial q}(q_i)\right)\cdot\delta q_i = 0.
\end{align*}
Since $\delta q_i$ is unconstrained for $1 \leq i \leq n-2$, 
\begin{align}
\frac{\partial T}{\partial \dot q}(\dot q_{i+1/2}) - \frac{\partial T}{\partial \dot q}(\dot q_{i-1/2}) = -h \frac{\partial V}{\partial q}(q_i), \quad i = 1, \ldots, n-2 \label{eq-DEL},
\end{align}
the \emph{discrete Euler-Langrange equations} of the system.

Unlike in the continuous settings, the discrete Euler-Lagrange equations do not always have a unique solution given initial values $q_0$ and $q_{1}$. Therefore in all that follows it is assumed that $T_d$ and $V_d$ are of a form so that (\ref{eq-DEL}) gives a unique $q_{i+1}$ given $q_i$ and $q_{i-1}$---this assumption always holds, for instance, in the typical case where $T_d$ is quadratic in $\dot q$. Then the discrete Euler-Lagrange equations give a well-defined discrete flow
\begin{align*}
F: (q_{i-1}, q_i) \mapsto (q_i, q_{i+1}),
\end{align*}
which recovers the entire trajectory from initial conditions, in perfect analogy to the continuous setting.

\subsection{Symplecticity of the Discrete Flow}

By analogy to the continuous setting, it is desired that $F$ preserve a symplectic form, just as $d\theta_L$ is preserved by $\Theta$. As in the continuous setting, trajectories $\q$ are restricted to those that satisfy the discrete Euler-Lagrange equations, and variations to first variations
(and the condition that these variations vanish at the endpoints is lifted), yielding
\begin{align*}
dS_d(\q) \cdot \delta \q = T'_d(q_{n-2}, q_{n-1})\cdot\delta q_{n-1} - T'_d(q_0,q_1)\cdot \delta q_0 - h \frac{\partial V}{\partial q}(q_{n-1}) \cdot \delta q_{n-1}.
\end{align*}
$F^k$ denotes the discrete flow $F$ composed with itself $k$ times, or $k$ ``steps'' of $F$. Again, all $\q$ satisfying (\ref{eq-DEL}) can be parametrized by initial conditions $\nu_0 = (q_0, q_1)$, and first variations by $\delta \nu_0 = (\delta q_0, \delta q_1)$, so that the discrete action can be rewritten as
\begin{align*}
S_{id}(\nu_0) = \sum_{i=0}^{n-2} h L_d( F^i\nu_0 ).
\end{align*}

Putting together all of the pieces,
\begin{align*}
dS_{id}(\nu_0) \cdot \delta\nu_0 &= dS_d(\q) \cdot \delta \q\\
&= T'_d(q_{n-2},q_{n-1})\cdot \delta q_{n-1} - T'_d(q_0,q_1)\cdot \delta q_0 - h \frac{\partial V}{\partial q}(q_{n-1}) \cdot \delta q_{n-1}\\
&= \left(T'_d(q_a,q_b) - h \frac{\partial V}{\partial q}(q_b)\right) d q_b \cdot (\delta q_{n-2},\delta q_{n-1})\Big\vert_{q_a=q_{n-2},\ q_b=q_{n-1}}\\
&= \left[T'_d(F^{n-2} \nu_0) - h V'(F^{n-2} \nu_0) \right] d q_b \cdot {F^{n-2}}_* \delta \nu_0 - T'_d(\nu_0)dq_a \cdot \delta \nu_0\\
&= \theta^+_{F^{n-2}\nu_0} \cdot {F^{n-2}}_* \delta \nu_0 + \theta^-_{\nu_0} \cdot \delta\nu_0\\
&= \left({F^{n-2}}^* \theta^+\right)_{\nu_0} \cdot \delta\nu_0 + \theta^-_{\nu_0} \cdot \delta\nu_0.
\end{align*}
for the indicated two-forms $\theta^+$ and $\theta^-$. Since $d(h L_d) = \theta^+ + \theta^-$, $d^2(h L_d) = 0 = d\theta^+ + d\theta^-$. Moreover the initial conditions $\nu_0$ are arbitrary, hence
\begin{align*}
d^2S_{id} = 0 = {\left(F^{n-2}\right)}^* d\theta^+ + d\theta^- = -{\left(F^{n-2}\right)}^*d\theta^- + d\theta^-,
\end{align*}
so
\begin{align*}
d\theta^- = {\left(F^{n-2}\right)}^* d\theta^-.
\end{align*}
Since $n$ is arbitrary, the discrete flow $F$ preserves the symplectic form $d\theta^-$. Using backwards error analysis, it can be shown that this geometric property guarantees that integrating with $F$ introduces no energy drift for a number of steps exponential in $h$~\cite{bHairer2006}, a highly desirable property when simulating molecular dynamic or other Hamiltonian systems whose qualitative behavior is substantially affected by errors in energy.

\section{Asynchronous Variational Integrators \label{s-AVI} }

In Section \ref{s-discrete} an action functional 
(\ref{eq-daction}) was formulated as the integration of a single discrete Lagrangian 
over a single time step size $h$. Such a construction is cumbersome 
when modeling multiple potentials of varying stiffnesses acting on 
different parts of the system: to prevent instability one must 
integrate the entire system at the resolution of the stiffest 
force. Asynchronous variational integrators (AVIs), introduced by Lew \etal~\cite{aLew2003}, are a family of numerical integrators, derived from a discrete Hamilton's principle, that support integrating potentials at different time steps. Their formulation assumes a spatial partition, with each potential depending only on the configuration of a single element; in this exposition, the general arguments set forth by Lew \etal\,are followed, but the notation and derivation departs from their work as necessary to support potentials with arbitrary, possibly non-disjoint spatial stencil.

Let $\{V^i\}$ be potentials with time steps $h^i$. Each potential $V^i$ is concerned with certain moments in time---namely, integer multiples of $h^i$---and these moments are inconsistent across triangles. Time is therefore subdivided in a way compatible with all triangles: for a $\tau$-length interval of time, the set $\Xi(\tau)$ is defined by
\begin{align*}
\Xi(\tau) = \bigcup_{V^i} \bigcup_{j=0}^{\lfloor \tau/h^i \rfloor} jh^i.
\end{align*}
That is, $\Xi(\tau)$ is the set of all integer multiples less than $\tau$ of all time steps. $\Xi$ can be ordered, and in particular let $\xi(i)$ be the $(i+1)$-st least element of $\Xi$. For ease of notation, also let $\omega^i(j) = \xi^{-1}(jh^i)$; that is, $\omega$ converts the $j$th timestep of potential $i$ into a global time.

If $n$ is the cardinality of $\Xi$, a trajectory of duration $\tau$ is then discretized by linearly interpolating intermediate configurations $q_0, q_1, \ldots, q_{n-1}$, where $q_i$ is the configuration of the system at time $\xi(i)$. Velocity is discretized as $\dot q_{k+1/2} = \frac{q_{k+1}-q_k}{\xi(k+1)-\xi(k)}$ on the segment of the trajectory between $q_k$ and $q_{k+1}$. A global action functional of these trajectories is needed, and can be constructed in the natural way:
\begin{align*}
S_{g}(\q) = \sum_{j=0}^{n-2}\left[\xi(j+1)-\xi(j)\right]T_d\left[q_j, q_{j+1}, \xi(j), \xi(j+1)\right]-\sum_{V^i}\sum_{j=1}^{\lfloor \tau/h^i\rfloor} h^i V^i (q_{\omega^i(j)}), \label{eq-actionAVI}
\end{align*}
where, for $T(\dot{q})$ the kinetic energy of the entire configuration, $T_d(q_a, q_b, t_a, t_b) = T\left(\frac{q_b-q_a}{t_b-t_a}\right)$. For use in the following, also let $T'_d(q_a, q_b, t_a, t_b) = \frac{\partial T}{\partial \dot{q}}\left(\frac{q_b-q_a}{t_b-t_a}\right)$.

No attempt has been made to define a Lagrangian pairing the kinetic and potential energy terms; it will be seen that an action defined in this way still leads to a multisymplectic numeric integrator. To this end, Hamilton's principle $dS_g(\q)\cdot \delta \q=0$ is imposed for variations $\delta \q = \{ \delta q_0, \ldots, \delta q_{n-1} \}$ with $\delta q_0 = \delta q_{n-1} = 0$. Then $S_g$ can be rewritten as
\begin{align}
S_{g}(\q) = \sum_{j=0}^{n-2}\left[\xi(j+1)-\xi(j)\right]T_d\left[q_j, q_{j+1}, \xi(j), \xi(j+1)\right]-\sum_{j=1}^{n-1}\sum_{h^i \vert \xi(j)}  h^i V^i (q_j),
\end{align}
where the notation $h^i \vert \xi(j)$ is abused to mean ``all indices $i$ for which $h^i$ evenly divides $\xi(j),$'' so that
\begin{align*}
dS_g(\q) \cdot \delta \q &= \sum_{j=0}^{n-2} T'_d\left[q_j, q_{j+1}, \xi(j), \xi(j+1)\right]\cdot \left(\delta q_{j+1} - \delta q_j\right) - \sum_{j=1}^{n-1} \sum_{h^i \vert \xi(j)} h^i \frac{\partial V_i}{\partial q}(q_j)\cdot \delta q_j\\
&= T'_d\left[q_{n-2}, q_{n-1}, \xi(n-2), \xi(n-1)\right]\cdot \delta q_{n-1} - T'_d\left[q_0, q_1, \xi(0), \xi(1)\right]\cdot \delta q_0\\
&\quad - \sum_{h^i \vert \xi(n-1)}h^i \frac{\partial V^i}{\partial q}(q_{n-1}) \cdot \delta q_{n-1}\\
&\quad + \sum_{j=1}^{n-2} \left( T'_d\left[q_{j-1}, q_{j}, \xi(j-1), \xi(j)\right] - T'_d\left[q_j, q_{j+1}, \xi(j), \xi(j+1)\right]-\sum_{h^i \vert \xi(j)} h^i \frac{\partial V^i}{\partial q}(q_j)\right)\cdot \delta q_j\\
&= \sum_{j=1}^{n-2} \left( T'_d\left[q_{j-1}, q_{j}, \xi(j-1), \xi(j)\right] - T'_d\left[q_j, q_{j+1}, \xi(j), \xi(j+1)\right]-\sum_{h^i \vert \xi(j)} h^i \frac{\partial V^i}{\partial q}(q_j)\right)\cdot \delta q_j.
\end{align*}

The Euler-Lagrange equations are then
\begin{align}
\frac{\partial T}{\partial \dot q}(\dot q_{k+1/2}) - \frac{\partial T}{\partial \dot q}(\dot q_{k-1/2}) = - \sum_{h^i \vert \xi(k)} h^i \frac{\partial V^i}{\partial q^i}(q_{k}), \label{eq-AVIDEL}
\end{align}
These equations are similar to those derived for synchronous variational integrators (\ref{eq-DEL}), except that only a subset of potentials $V_d^i$ contribute during each time step. As in the synchronous case, if, as is typical, $T_d(\dot{q})$ is quadratic in $\dot q$, the system (\ref{eq-AVIDEL}) gives rise to an explicit numerical integrator that is particularly easy to implement in practice. Algorithm \ref{alg-AVI} gives pseudocode for such integration when $T_d = \dot q^T \M \dot q$ for a mass matrix $\M$; Lew \etal~\cite{tLew2003} discuss the algorithm in greater detail.

\begin{algorithm}
\caption{An algorithm for integrating the trajectory given by the AVI Euler-Lagrange equations (\ref{eq-AVIDEL}) adapted from Lew \etal~\cite{tLew2003}}
\label{alg-AVI}
\begin{algorithmic}
\STATE Let events be (potential, time step, time) triplets $E = (V, h, t)$.
\STATE Denote by $q_V$ the configuration subspace on which $V$ depends.
\STATE Let $PQ$ be a priority queue of events, sorted by event times $E.t$.
\STATE $T_g \leftarrow 0$ \COMMENT{$T_g$ maintains the value of the simulation clock}
\STATE $q \leftarrow q_0$ \COMMENT{Set up initial conditions}
\STATE $\dot q \leftarrow \dot q_0$
\FORALL{$V_i$}
\STATE $E_i \leftarrow (V_i, h^i, h^i)$ \COMMENT{Add all potentials to the queue as events}
\STATE $PQ$.push($E_i$)
\ENDFOR
\LOOP
\STATE $(V, h, t) \leftarrow PQ$.pop
\STATE $q \leftarrow q + (t-T_g) \dot q$ 
\STATE $\dot q_V \leftarrow \dot q_V - h M_V^{-1} \frac{\partial V}{\partial q_V}$ \COMMENT{Update only those elements affected by this event.}
\STATE $PQ$.push($V, h, t + h$) \COMMENT{Return the event to the queue, with a new, later time}
\STATE $T_g \leftarrow t$ \COMMENT{Update the simulation clock}
\ENDLOOP
\end{algorithmic}
\end{algorithm}

\subsection{Multisymplecticity} \label{s-AVImulti}

The right hand side of (\ref{eq-AVIDEL}) depends on $\xi(k)$, and so unlike (\ref{eq-DEL}), the Euler-Lagrange equations for AVIs are time dependent, and do not give rise to a uniform update rule $F(q_{i-1}, q_{i})\mapsto (q_{i}, q_{i+1})$. Instead, consider the total, time-dependent flow $\hat{F}^k(q_{0}, q_{i}) \mapsto (q_{k}, q_{k+1})$. Once again, trajectories satisfying (\ref{eq-AVIDEL}) are parametrized by $\nu_0=(q_0, q_1)$, and first variations by $\delta \nu_0=(\delta q_0, \delta q_1)$. By restricting to such trajectories and variations, the action (\ref{eq-actionAVI}) can be rewritten as
\begin{align*}
S_{\textrm{iAVI}}=\sum_{j=0}^{n-2} \left[\xi(j+1)-\xi(j)\right]T_d\left(\hat{F}^j(\nu_0), \xi(j), \xi(j+1)\right) - \sum_{V^i} \sum_{j=0}^{\lfloor \tau/h^i \rfloor-1} h^i V^i_d (\hat{F}^{\omega^i(j+1)}(\nu_0)).
\end{align*}

Then, for ${V^i_d}'(q_a,q_b) = \frac{\partial V^i}{\partial q}(q)$, 
\begin{align*}
dS_{\textrm{iAVI}}(\nu) \cdot \delta \nu &= dS_{g}(\q) \cdot \delta \q \\
&= T'_d\left[q_{n-2},q_{n-1},\xi(n-2),\xi(n-1)\right]\cdot \delta q_{n-1} - T'_d\left[q_0,q_1,\xi(0),\xi(1)\right]\cdot \delta q_{0}\\
&\quad - \sum_{V^i}\sum_{h^i \vert \xi(n-1)} h^i \frac{\partial V^i}{\partial q^i}(q_{n-1})\cdot \delta q_{n-1}\\
&= T'_d\left[\hat{F}^{n-2}(\nu_0), \xi(n-2), \xi(n-1)\right] \cdot \delta q_{n-1} - T'_d\left[\nu_0, \xi(0), \xi(1)\right]\cdot \delta q_{0}\\
&\quad - \sum_{V^i} \sum_{h^i \vert \xi(n-1)} h^i {V^i_d}'\left[\hat{F}^{n-2}(\nu_0)\right]\cdot \delta q_{n-1}\\
&= \theta^-_{\nu_0} \cdot \delta \nu_0 + \theta^+_{\hat{F}^{n-2}\nu_0} \cdot {\hat{F}^{n-2}}{}_* \delta \nu_0\\
&= (\theta^- + {\hat{F}^{n-2}}{}^* \theta^+)_{\nu_0} \cdot \delta \nu_0
\end{align*}
for one-forms $\theta^-$ and $\theta^+$. Once again
\begin{align}
0 = d^2 S_{\textrm{iAVI}} = d\theta^- + \hat{F}^{n-2}{}^*d\theta^+, \label{eq-multi}
\end{align}
but unlike when the action was a sum of Lagrangians, from the \emph{multisymplectic form formula} (\ref{eq-multi}) there is no way of relating $d\theta^-$ to $d\theta^+$, and thus discrete symplectic structure preservation is not recovered. Nevertheless, Lew \etal~\cite{aLew2003} conjecture that this multisymplectic structure leads to the good energy behavior observed for AVIs.

\section{Discrete Penalty Layers}

The above reformulation of AVIs can be leveraged to resolve collisions with guaranteed perfect robustness, and via momentum-symplectic integration, so that the energy behavior of the system as a whole remains good. Consider a standard penalty force approach, which for every two elements $A, B$ and surface thickness $\eta$ defines the gap function
\begin{align*}
g_{\eta}(q) = \inf_{a\in A, b\in B} \|a-b\|-2\eta
\end{align*}
measuring the proximity of $A$ to $B$. 

\begin{figure}[t]
\centering
\includegraphics[width=3in]{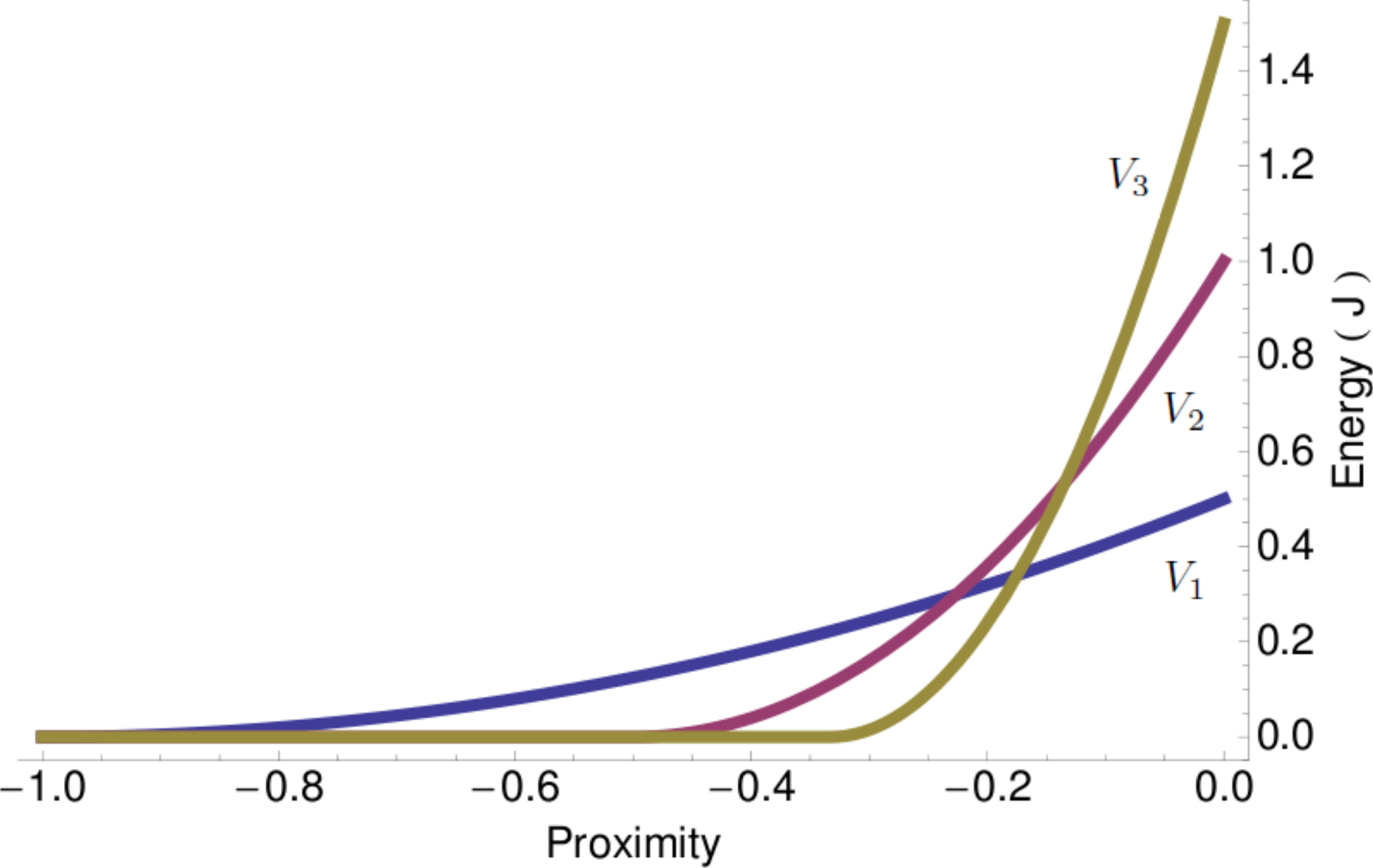}
\hfill
\includegraphics[width=3in]{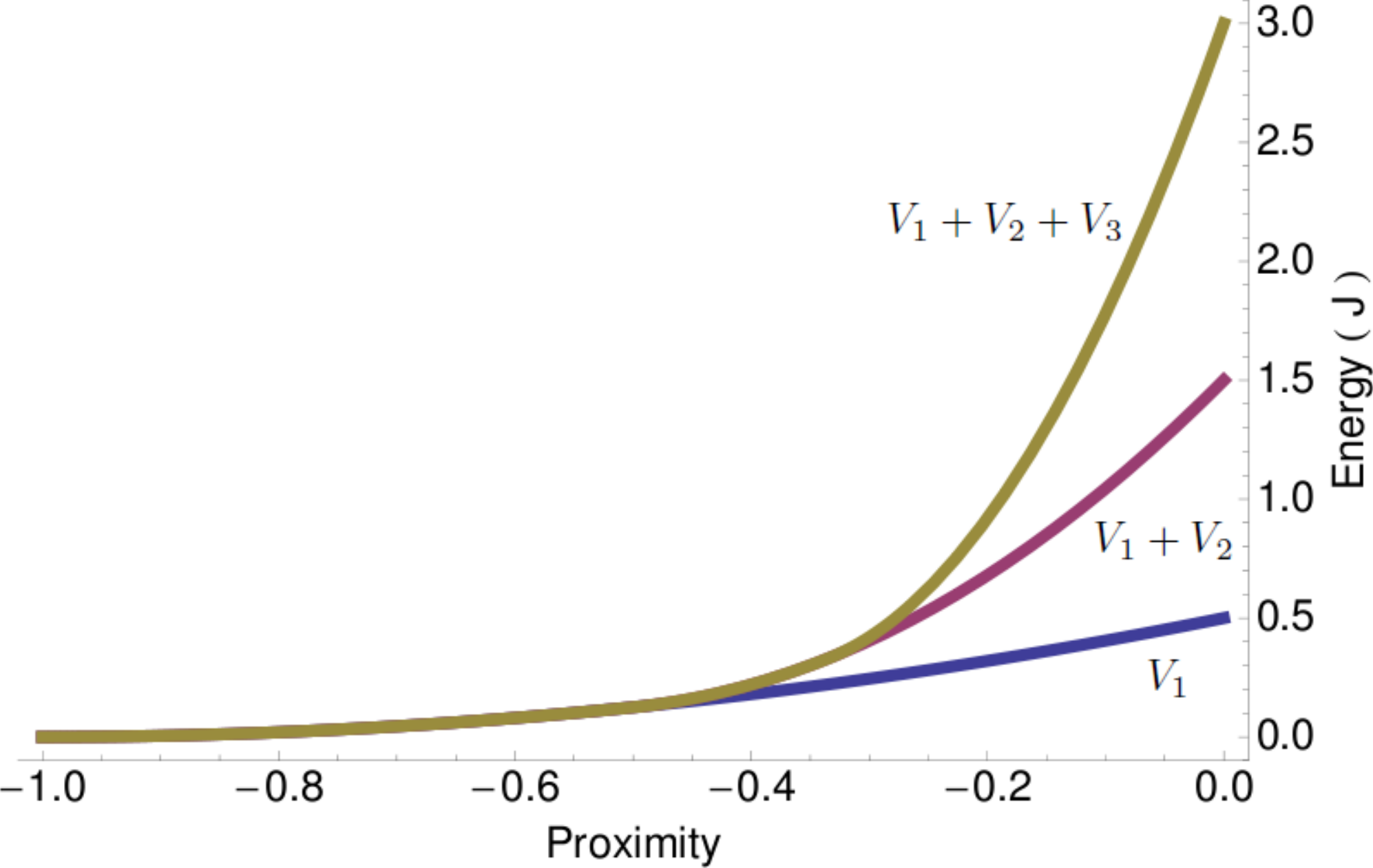}
\caption{Plots of the potential energy of the first three layers as a function of gap function $g$ (left), and a plot of the total potential energy contributed by all layers $\leq n$ for $n=1, 2, 3$ (right). Notice the potential energy diverges as separation distance approaches $0$, guaranteeing that collision response is robust.}
\label{fig-pots}
\end{figure}

The penalty potential is then defined as
\begin{align*}
V(q) =
\begin{cases}
0 & g_{\eta}(q) > 0\\
kg_{\eta}(q)^2 & g(q) \leq 0,
\end{cases}
\end{align*}
where $k$ is a user-specified stiffness. As previously discussed, $V$ alone does not robustly prevent interpenetrations: the potential can be viewed as placing a spring between the approaching elements, and for sufficiently large relative momentum in the normal direction, the spring will fully compress, then fail. However, consider placing an infinite family of potentials $V_l$, $l = 1, 2, \ldots,$ between the primitives, where
\begin{align*}
V_l(q) =
\begin{cases}
0 & g_{\eta/l}(q) > 0 \\
l^3 k g_{\eta/l}^2 & g_{\eta/l}(q) \leq 0.
\end{cases}
\end{align*}
The region $\frac{\eta}{n+1} \leq d(q) \leq \frac{\eta}{n}$, where exactly $n$ of the potentials are active, is called the $n$-th \emph{discrete penalty layer}. Figure \ref{fig-pots} shows a plot of the potential energy of the first few potentials for the case $\eta=k=1$, as well as the cumulative potential energy of all of the potentials.

The total potential energy of the springs when fully compressed is
\begin{align*}
\sum_{l=1}^\infty l^3 k 4 \left(\frac{\eta}{l}\right)^2 = 4 k\eta^2\sum_{l=1}^\infty l,
\end{align*}
which diverges. The infinite array of potentials is guaranteed to stop all collisions. This guarantee in no way depends on the chosen stiffness $k$: although performance and error will vary with the choice of stiffness, unlike for penalty forces the stiffness does not affect the guarantee. The method is always guaranteed to be robust. 

There is one obstruction to implementing this scheme in practice: integrating the $l$-th spring stably and with good energy behavior requires a time step proportional to $\frac{1}{l^{3/2}}$, which vanishes as $l\to\infty.$ Using a traditional integrator, one could decide ahead of time to only simulate the first few springs---but then the guarantee that no penetrations will occur is lost, and the simulation must be run at a prohibitively small time step. AVIs, with the above modifications, and a bit of extra bookkeeping, are a first step towards alleviating the problem, by allowing the user to assign each spring its own time step. This bookkeeping is now described, in terms of modifications to the basic Algorithm \ref{alg-AVI}.

\section{The Asynchronous Algorithm}
\label{sKDS}
AVIs allow each penalty layer to be assigned a different time step, so that less stiff ($l$ small) layers can take large time steps regardless of the presence of the stiffer layers. However, it is still not possible as a practical matter to integrate the system, since arbitrarily large $l$ would need arbitrarily small time steps, and the global time in Algorithm \ref{alg-AVI} would never advance. The following observation surmounts this obstacle: at any time during a well-posed simulation, the number of layers that are exerting a non-zero force, or that are \emph{active}, is finite. More precisely, a simulation is well-posed if its total energy over time is bounded---that is, if the simulation begins in a non-penetrating state; all prescribed, infinite-mass bodies are stationary; and only a finite amount of energy is added over time in the form of external forcing. Inactive penalty potentials can be ignored by Algorithm \ref{alg-AVI} entirely, since they do not change configurational velocity, and the position integration that would take place during the handling of an inactive potential can just as well be done by the following event. Therefore the simulation would be guaranteed to never stop making progress if there is a lower bound for the amount of global time $T_g$ that elapses with the processing of any event. Such a lower bound exists if there is a way to detect which penalty potentials are active or inactive at all times and remove all inactive events from the priority queue $PQ$.

Suppose that at the start of the simulation, all penalty layers are inactive. Thus no penalty layer events are needed on the queue. For each pair of simulation elements, the time $t_a$ that the first penalty layer would become active (assuming all elements continue along the trajectory described by their initial velocities) can be calculated, and the corresponding event added to the queue at that time. Such an approach suffers from two problems, however. Firstly, solving for the time when the gap function will be zero is easy in some cases, such as if the elements are two spheres or two planes, but can involve expensive root solves in others, such as if the elements are two non-rigid triangular elements of a thin shell simulation. Secondly, the times computed are fragile: should any event alter the velocity of one of the elements (such as a material force, or gravity, or another penalty force if one of the elements collides with a third party) the activation time is no longer valid and must be recomputed.

Instead of an exact time, only a conservative guarantee, or \emph{certificate}~\cite{aGuibas1998}, that the first penalty layer will not be active before some time $t_c$ (where necessarily $t_c \leq t_a$) is truly needed. For example, one certificate is the existence of an $2\eta$-thick \emph{planar slab} $S$ that separates the two elements up until time $t_c$, where $\eta$ is the thickness of the first penalty layer. For an $m$-dimensional configuration space, such a planar slab is understood to be an extrusion of an $(m-1)$---dimensional affine subspace. Concretely, let $w$ be a unit vector in $\mathbb{R}^m$, $w_i$ be $m-1$ linearly independent vectors in $\mathbb{R}^{m}$ orthogonal to $w$, and $p$ a point in $\mathbb{R}^m$. Then the slab $S_{w,p}$ is the set
\begin{align*}
S_{w,p} = \Big\{ p+\alpha w + \sum_i \beta_i w_i \Big\vert -\eta \leq \alpha \leq \eta, \beta_i \in \mathbb{R}\Big\}.
\end{align*}
If such a slab separates the two elements, the first penalty layer cannot become active before $t_c$. This certificate can be placed as an event on the queue, with time $t_c$. The certificate might then suffer several fates:~\cite{aHarmon2009}
\begin{itemize}
\item{An event modifies the velocity of one of the elements before time $t_c$. The certificate placed on the queue is then no longer valid until time $t_c$, but instead until a new time $t_c'$ which may be sooner or later than $t_c$. The algorithm must thus \emph{reschedule} the certificate, by removing its event from the queue, and reinserting it at the appropriate new time.}
\item{The certificate event is popped from the queue without incident, but it is possible and convenient to find a new separating slab that guarantees the penalty layer does not activate before time $t'_c > t_c$. This new certificate can then be pushed on the queue for time $t'_c$.}
\item{The certificate event is popped from the queue without incident, but finding a new slab is impossible, costly, or a slab can be found, but the new time $t'_c$ is judged heuristically to be too near $t_c$. The first penalty layer may then be activated early: doing so affects the efficiency, but not the correctness, of the simulation. Simultaneously, the algorithm searches for an $\eta$-thick separating slab to serve as a certificate that layer two is not yet active, and the whole process described above is repeated.}
\end{itemize}

Detecting when a penalty layer event becomes inactive, and should be removed from the queue, is much simpler than detecting layer activation: whenever a penalty force for layer $n$ is integrated, the algorithm simply checks if the force applied was 0. If so, and if the two elements in question are separating, layer $n$ is now inactive: it is not pushed back onto the queue (and instead a separating slab of thickness $\eta/n$ is sought.)

It is very important to note that when an event becomes active and is added back into the event priority queue, it is done so at a time that is \emph{an integer multiple of its timestep from its last time of integration.} That is, those times when integration would do nothing have been optimized away, but the potential's ``integration clock'' has not been tampered with or realigned, since every potential having a fixed-size time step was fundamental to the proof that asynchronous variational integration is multisymplectic. The spring-on-a-plane example described below underlines the danger of failing to maintain such a fixed time step.

For an event $E$, denote all simulation elements on which $E$ depends the \emph{support} of $V$. Denote all simulation elements whose velocities are modified by $E$ the \emph{stencil} of $E$. For force integration events, there is no distinction between stencil and support. Certificates have a support, but no stencil. Algorithm \ref{alg-AVIslab} uses this terminology to incorporate the above into the AVI algorithm.

\begin{algorithm}
\caption{Proposed algorithm for asynchronous contact resolution.}
\label{alg-AVIslab}
\begin{algorithmic}
\STATE Let force events be (potential, time step, time) triplets $E = (V, h, t)$.
\STATE Let $PQ$ be a priority queue of events, sorted by event times $E.t$.
\STATE $T_g \leftarrow 0$ \COMMENT{$T_g$ maintains the value of the simulation clock}
\STATE $q \leftarrow q_0$ \COMMENT{Set up initial conditions}
\STATE $\dot q \leftarrow \dot q_0$
\STATE Push non-penalty (e.g. material) events on the queue
\FORALL {pairs of elements $K_1$, $K_2$}
\STATE $E \leftarrow$ FindCertificate($K_1$, $K_2$)
\STATE $PQ$.push($E$)
\ENDFOR
\LOOP
\STATE $E \leftarrow PQ$.pop
\STATE $q \leftarrow q + (E.t-T_g) \dot q$
\IF{$E$ is a force event}
\STATE handleForceEvent($PQ$, $E$)
\ELSE 
\STATE handleCertificateEvent($PQ$, $E$)
\ENDIF
\STATE $T_g \leftarrow E.t$ \COMMENT{Update the simulation clock}
\ENDLOOP
\end{algorithmic}
\end{algorithm}

\begin{algorithm}
\caption{handleForceEvent}
\begin{algorithmic}
\REQUIRE Priority queue of events $PQ$ and force event $E$ that needs processing
\STATE \COMMENT{Processing a force event $E$ is a three-step process: integrating the force, rescheduling all events whose support depends on $E$'s stencil, and lastly, resceduling $E$ itself.}
\FORALL{$i$ in Stencil($E$)}
\STATE $\dot q_i \leftarrow \dot q_i - (E.h) M_i^{-1} \frac{\partial E.V}{\partial q_i}$ \COMMENT{Update only those elements affected by this event.}
\ENDFOR
\STATE \COMMENT{Reschedule all events whose support depends on $E$'s stencil}
\FORALL{Certificate events $E'$ with Stencil($E$) $\cap$ Support($E'$) $\neq \emptyset$}
\STATE{$PQ$.remove($E'$)}
\STATE Schedule($E'$)
\STATE $PQ$.push($E'$)
\ENDFOR
\STATE \COMMENT{If $E$ was a penalty force event, it exerted 0 force, and the two primitives in question are separating, then we no longer need it}
\IF{$E$ is a penalty force event and $\frac{\partial E.V}{\partial q_i} = 0$}
\IF{$E.V.K1$ and $E.V.K2$ have positive relative velocity (are separating)}
\RETURN
\ENDIF
\IF{addCertificate($E.V.K1, E.V.K2$)}
\RETURN
\ENDIF
\ENDIF
\STATE \COMMENT{Otherwise, reschedule $E$ itself}
\STATE $PQ$.push($V, h, t + h$)
\end{algorithmic}
\end{algorithm}

\begin{algorithm}
\caption{handleCertificateEvent}
\begin{algorithmic}
\REQUIRE Priority queue of events $PQ$ and certificate event $E$ that needs rescheduling
\IF{not addCertificate($E.K1$, $E.K2$)}
\STATE \COMMENT {Finding a new certificate failed. We must thus activate a penalty force, one layer deeper than the deepest currently active penalty force event.}
\STATE $CurLayer \leftarrow \max_{\{\textrm{penalty events $E'$ on queue for $E.K1$ and $E.K2$}\}} E'.layer$
\STATE $E' \leftarrow $ new PenaltyForceEvent($E.K1$, $E.K2$, $CurLayer+1$) 
\STATE $PQ$.push($E'$) \COMMENT{Push the appropriate penalty force event on the queue}
\ENDIF
\end{algorithmic}
\end{algorithm}

\begin{algorithm}
\caption{addCertificate}
\begin{algorithmic}
\REQUIRE Priority queue of events $PQ$, and two elements $K1$ and $K2$

\STATE \COMMENT{Attempts to find a certificate for the collision of $K1$ against $K2$ and add it to the queue. Returns true if one was found.}
\STATE $E' \leftarrow $ FindCertificate($K1$, $K2$)
\IF{$E'$ was successfully found}
\STATE $PQ$.push(FindCertificate($K1$, $K2$))
\RETURN \TRUE
\ENDIF
\RETURN \FALSE
\end{algorithmic}
\end{algorithm}

In Algorithm \ref{alg-AVIslab} and its accompanying subalgorithms, the behavior of the functions FindCertificate and Schedule will depend on the type of certificate chosen. FindCertificate returns a new certificate for a given pair of elements, if possible and practical, and Schedule computes the time a certificate becomes invalid, as described in the paragraphs above. For thin shell simulation, where all simulation elements are convex triangles, edges, and vertices, separating slabs serve as ideal certificates, since it is cheap to compute Schedule, in this case by calculating element-plane intersection times. Although any choice of certificate, and heuristic for when to abort searching for a new certificate, preserves the correctness of the algorithm, the progress property described in the first paragraph of this section relies on the certificates efficiently weeding out inactive events so that some certificate is found before all (infinitely many) layers for a pair of elements are activated. No problems have been observed using separating slabs for thin-shell simulations, but different certificates may be needed, e.g., for concave rigid bodies.

\subsection{Further Optimizations}
The technique explored in the previous section, of finding a sequence of conservative certificates guaranteeing that some property holds, instead of calculating an exact time when that property stops holding, is the central idea behind a wide class of algorithms known as \emph{Kinetic Data Structures} (KDSs) \cite{aGuibas1998}. In the case described above, the property was inactivity of a given penalty layer. KDSs are particularly well-suited for an asynchronous approach, since certificate expiration times may not all align to some convenient simulation clock, and the required rescheduling of certificates/searching for new certificates can reuse the priority queue data structure already needed for force integration events. To improve the efficiency of the implementation used to create the examples below, several more KDSs in addition to the separating slabs discussed above were implemented: a bounding volume hierarchy \cite{aKlosowski1998} was used to take advantage of the fact that spatially distant elements are unlikely to collide, separation lists \cite{aWeller2006} to optimize the bookkeeping of this hierarchy, and a novel KDS was devised to leverage the observation that high-frequency, low amplitude oscillations in velocity do not significantly change a separation slab's expiration time, so that rescheduling is in many cases unnecessary. All of the improvements are described in greater detail in \cite{aHarmon2009}.

\section{Dissipation}

The framework, as described so far, gives near-perfect long-time energy conservation. In the real world, however, many dissipative phenomena are observed --- for instance friction, spring damping, and non-unit coefficients of restitution during collisions. 
Several simple modifications can be made to the proposed method to take such dissipation into account.
Qualitatively, these have performed well in practice: energy seems to behave well over long times for dissipative systems analogously to the near-conservation observed for Hamiltonian systems, but a theoretical understanding of this good behavior remains future work.

\subsection{Coefficient of Restitution}
It is often desirable to simulate semi-elastic or inelastic collisions. A simple modification to the potential $V_l$ allows the use of arbitrary coefficients of restitution $e$:
\begin{align*}
V_l(q) =
\begin{cases}
0 & g_{\eta/l}(q) > 0\\
l^3 ks g_{\eta/l}(q)^2 & g_{\eta/l}(q) \leq 0,
\end{cases}
\end{align*}
where $s$ is $e$ if the primitives are separating, $1$ otherwise. The penalty layers exert their full force during compression, then
weaken according to the coefficient of restitution during decompression. 

This approach, while simple, does have a limitation in the inelastic limit $e=0$: due to error introduced by numerical integration, two colliding primitives may have non-zero, though small, post-response normal relative velocity. The magnitude of this velocity is at most $k \eta/l $, so it can be limited by choosing a small enough base stiffness $k$.

\subsection{Friction}

The Coulomb friction model is a simple approximation to kinetic friction: at a point of contact between two bodies, the Coulomb force has magnitude $\mu | F_n |$, where $\mu$ is a coefficient of friction and $F_n$ is the normal force at the contact points, and has direction opposite the relative tangential motion of the contact points.

Whenever an impulse is applied during integration of a penalty layer, a corresponding frictional impulse can also be applied. Just as increasingly stiff penalty forces are applied
for contact forces, friction forces are increasingly applied (equal
to $\mu |F_n|$) to correctly halt high-speed tangential motion. Notice that these friction forces, like the material and contact penalty forces, are applied asynchronously: every layer applies friction independently at its own time step.

This simple, asynchronous formulation of friction fits very naturally into the framework of AVIs. Unfortunately, it is unsuitable for simulations featuring static friction, such as a block of wood resting on an inclined plane. The above formulation, with friction applied piecemeal during penalty integration, is reactive instead of proactive, and in simulations of this type the block of wood has been observed ``creeping'' down the incline no matter how high a coefficient of restitution is chosen. A more comprehensive model of friction compatible with the AVI framework, which correctly handles static friction, remains future work.

\section{Results}

\subsection{Spring on a plane}

\begin{figure}[t]
\centering
\includegraphics[width=3in]{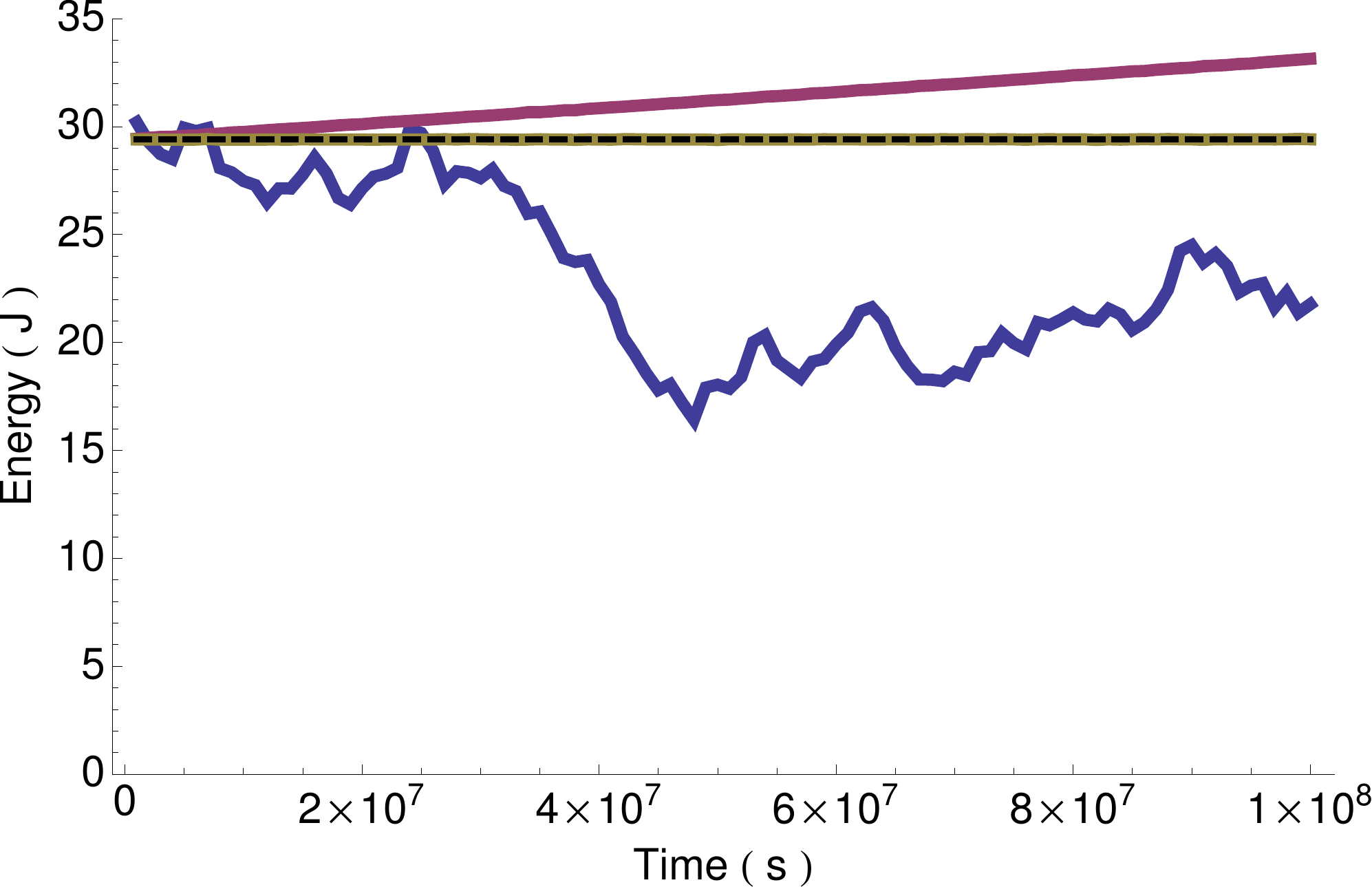}
\caption{The total energy of a spring bouncing on a plane integrated by mixing variational material force integration with impulse-based collision response (blue), the energy of this simulation when the alignment of the integration clocks is not respected (maroon), and the energy of this simulation when using the proposed method (brown). The exact energy is shown as a dashed black line.}
\label{fig-spring}
\end{figure}

As a simplest didactic demonstration of the proposed method, three experiments were conducted. A vertical spring of unit rest length, stiffness, and endpoint masses began each of the three simulations stationary a unit height above a fixed horizontal plane. The springs fell under a gravitational force of strength $1 m/s^2$, with impact handled in one of three different ways:

In the first experiment, gravity and the stretching force were integrated synchronously, and an instantaneous impulse was applied whenever the bottom of the spring touched the plane. Figure \ref{fig-spring} shows energy over time when using this method, compared to expected perfect energy conservation. Energy in this experiment, far from being conserved, took a random walk. 

In the second experiment, all forces were integrated asynchronously using the proposed method. The thickness $\eta$ was chosen to be $0.1$, and the penalty base stiffness $k$, $1000$. Energy in this case was well-conserved over long time: although energy experiences high-frequency, low-amplitude oscillations, there was no drift.

The importance of respecting the integrity of each potential's integration clock is highlighted in the third experiment. Instead of adding a force event onto the priority queue at an integer multiple of its time step, the events are added at the precise moment when each layer becomes active. As can be seen from the resulting energy plot (Figure \ref{fig-spring}), energy is no longer well-conserved, but instead seems to increase monotonically over time.

\subsection{Balls of Particles}
\begin{figure}[htb]
\centering
\includegraphics[height=2in]{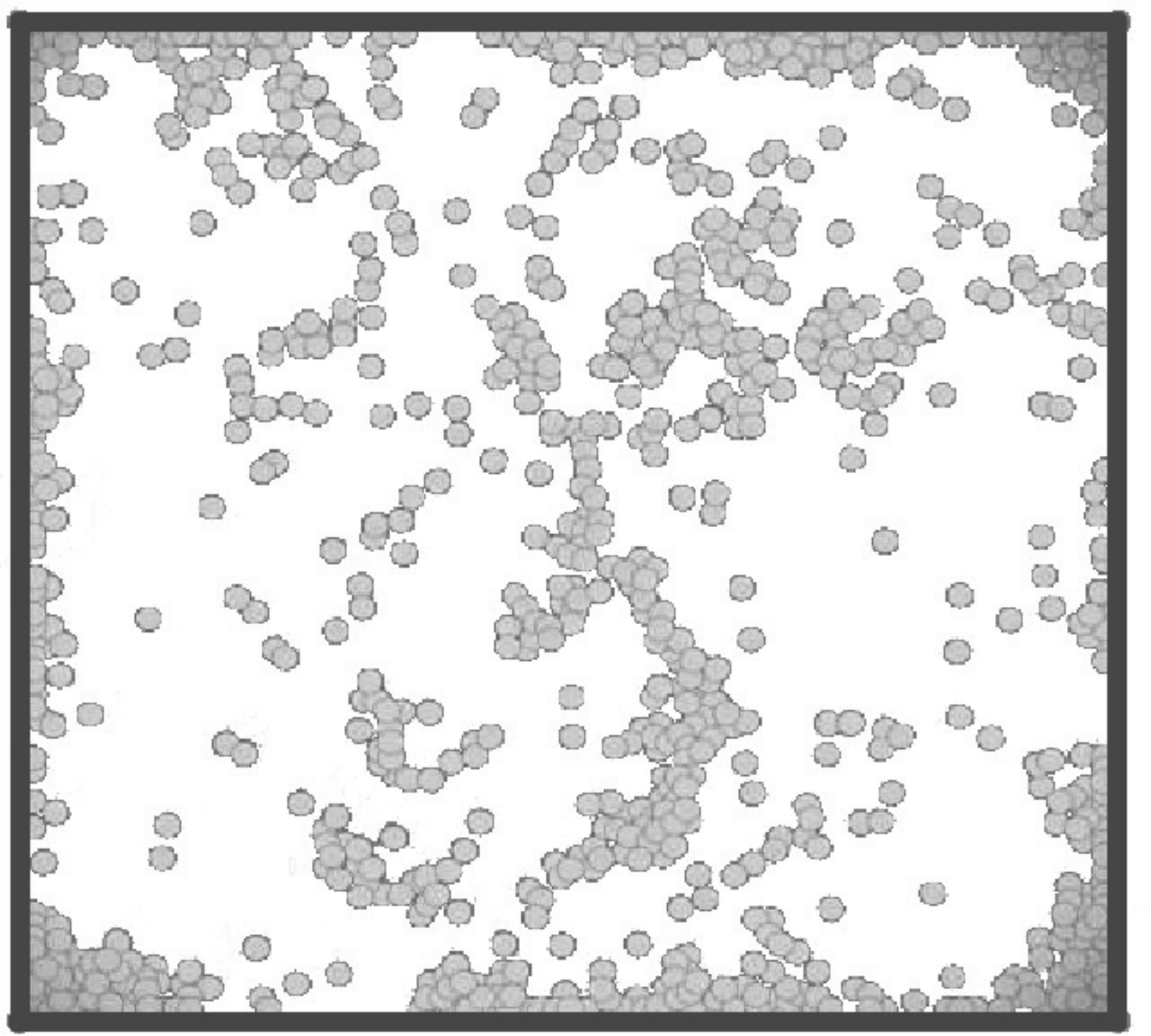}
\caption{A rigid box containing 900 spheres with random initial velocity, several minutes after the start of the simulation.}
\label{fig-box}
\end{figure}

\begin{figure}[htb]
\centering
\includegraphics[width=3in]{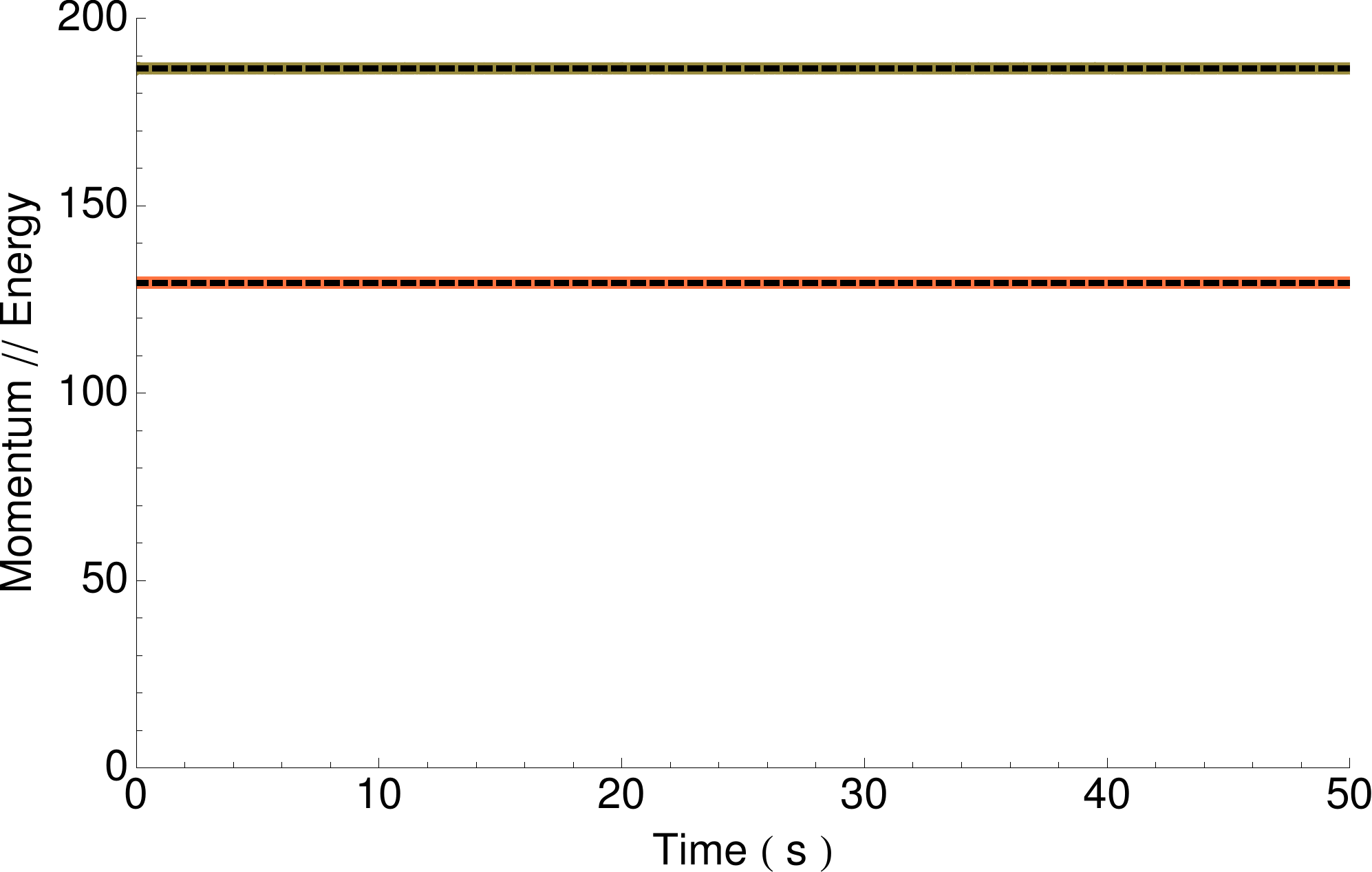}
\includegraphics[width=3in]{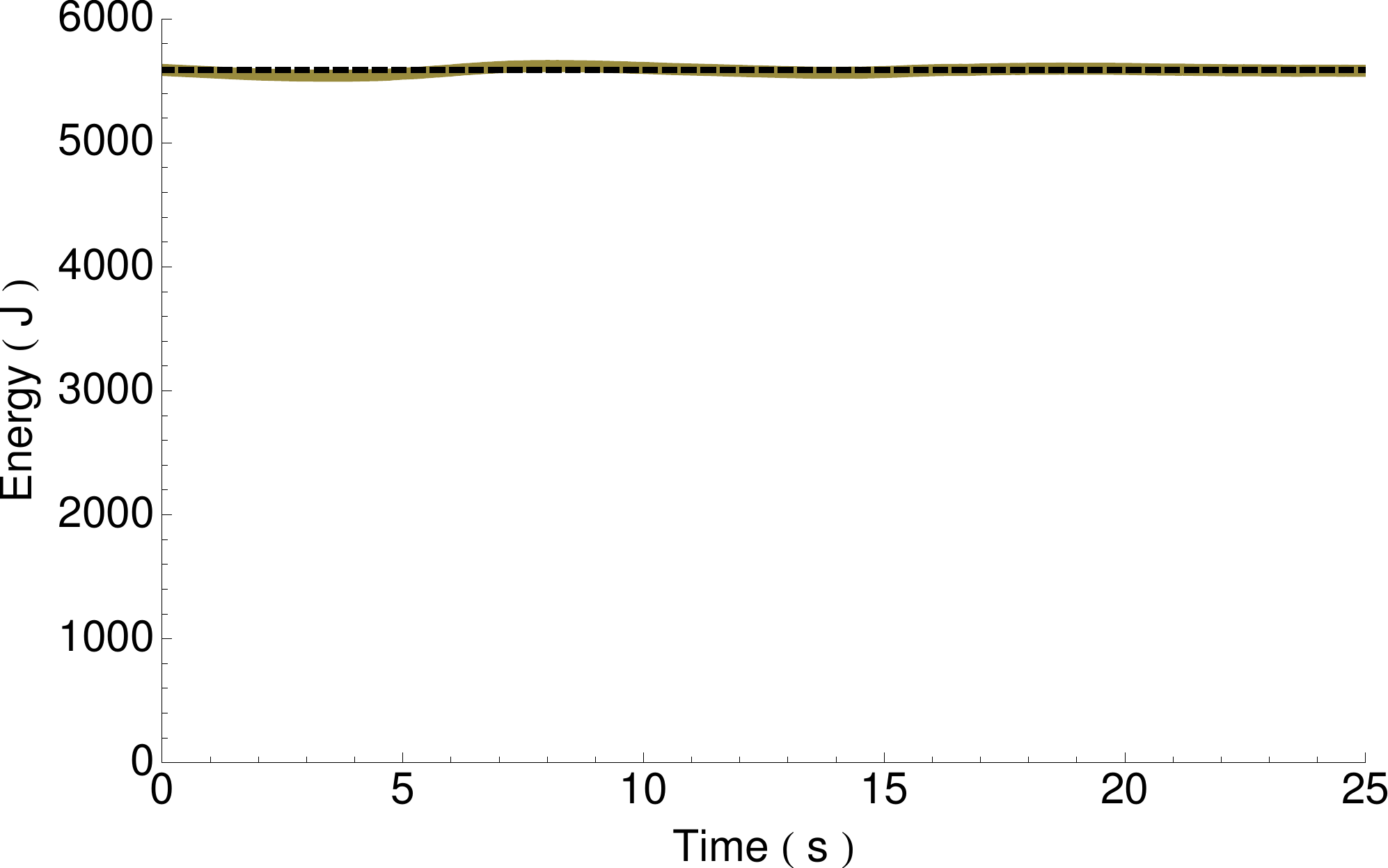}
\caption{Left: The energy (brown) and momentum (orange) of the box of the spheres as a function of time. Right: The energy of the box, with gravity added, over time.}
\label{fig-data}
\end{figure}

As an example that involves more collisions, consider a fixed $3 \unit{m} \times 3 \unit{m}$ square box. Inside this box $900$ spheres of radius $10 \unit{cm}$ were uniformly distributed, each of which was given a random velocity of magnitude between $0$ and $10 \unit{m/s}$. Figure \ref{fig-box} depicts this box after several minutes have elapsed. Energy over time is plotted in Figure \ref{fig-data} (left), and it is again almost perfectly conserved. Figure \ref{fig-data} (center) plots the magnitude of total momentum of the box over time, and it is exactly conserved, as expected since a multisymplectic integrator is used. Gravity ($9.8 \unit{m/s^2}$) was added to the box and again energy plotted over time (Figure \ref{fig-data}, right), and again good energy behavior was observed.

\begin{figure}[htb]
\centering
\includegraphics[width=3in]{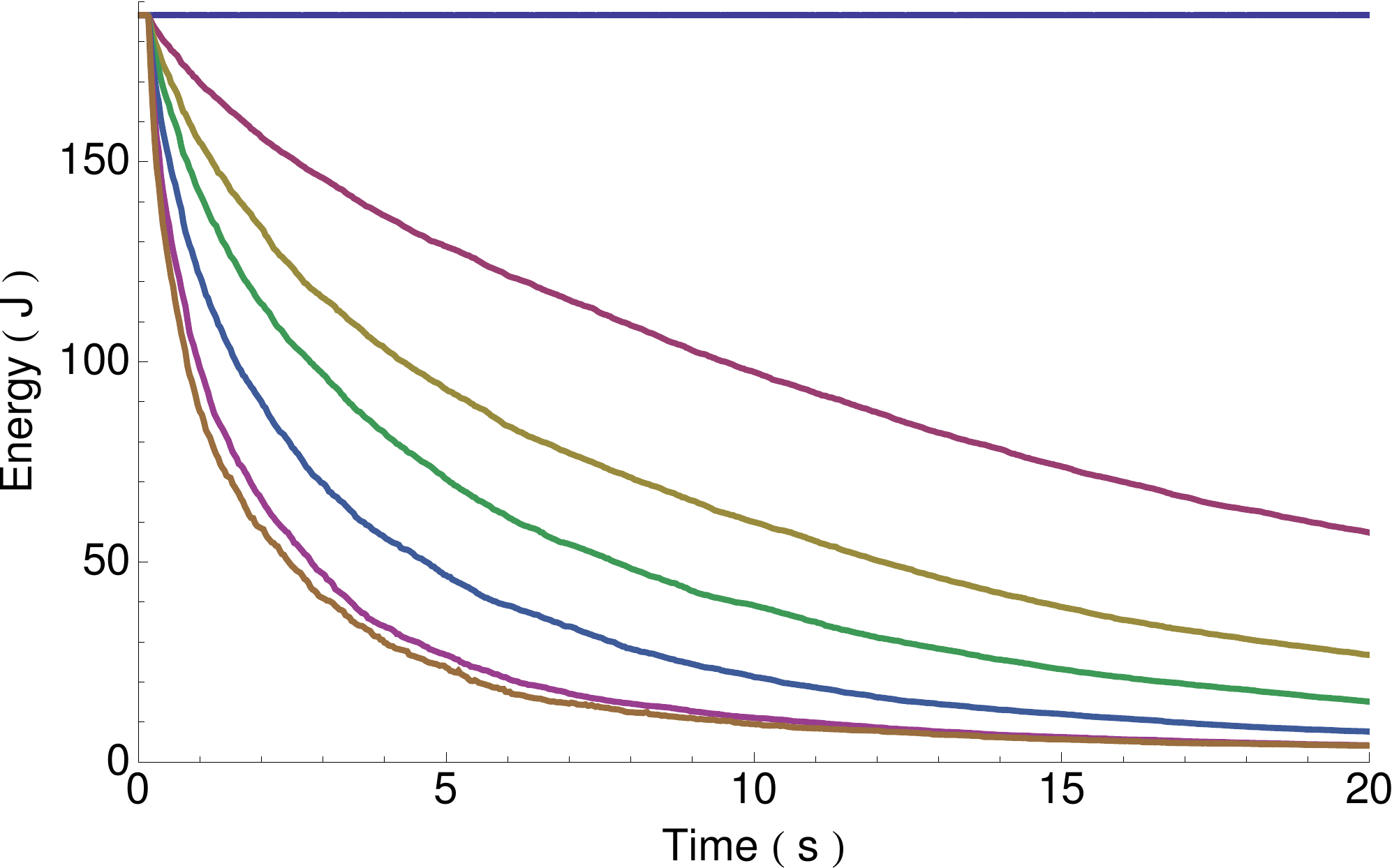}
\caption{The energy of the box of 900 spheres under different coefficients of restitution: from top to bottom, 1.0, 0.9, 0.8, 0.7, 0.5, 0.2, 0.0.}
\label{fig-cor}
\end{figure}

As a test of controllable dissipation by using a coefficient of restitution, the box with gravity was resimulated several times using different coefficients of restitution. Figure \ref{fig-cor} shows the resulting energy plots. For any chosen coefficient of restitution, the non-conservative energy behavior is qualitatively as one would expect.

\subsection{Sphere-Plate Impact}

The impact experiment of a spherical shell against a thin plate, as described in Cirak and West's article on Decomposition Contact Response (DCR)~\cite{aCirak2005}, was reproduced using the proposed framework. A sphere of radius $12.5 \unit{cm}$ approaches a plate of radius $35 \unit{cm}$ with relative velocity $100 \unit{m/s}$. Both the sphere and the plate have thickness $0.35 \unit{cm}$.
The time steps of the material forces (stretching and bending) are $10^{-7} \unit{s}$ (the same as those chosen by Cirak and West.)

\begin{figure}[htb]
\begin{center}
\includegraphics[width=3in]{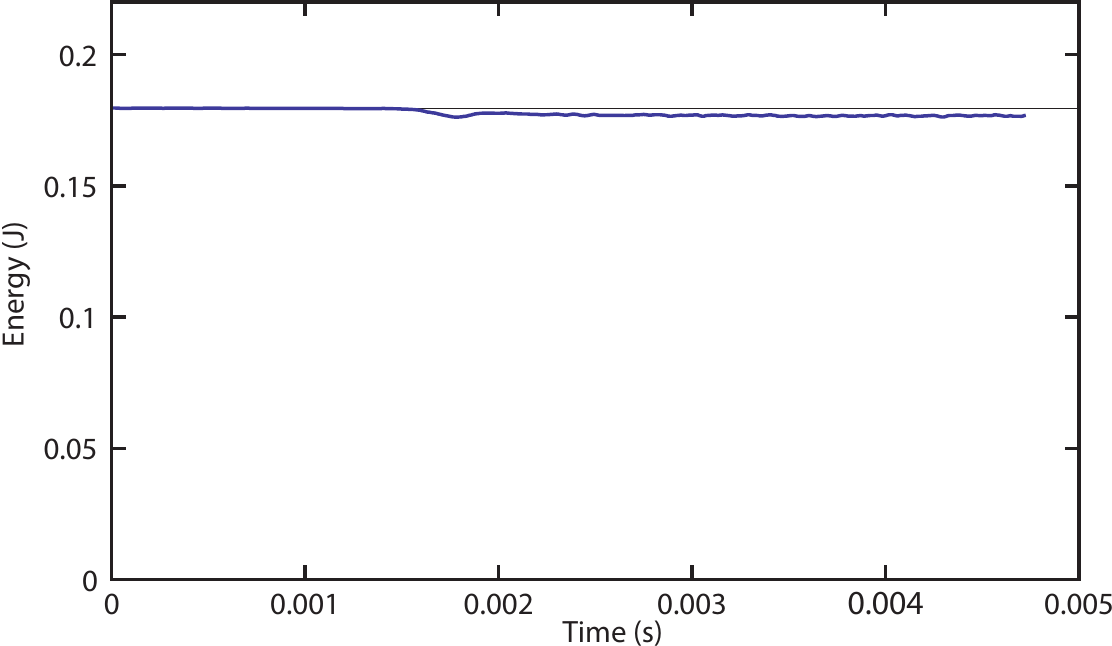}
\includegraphics[width=3in]{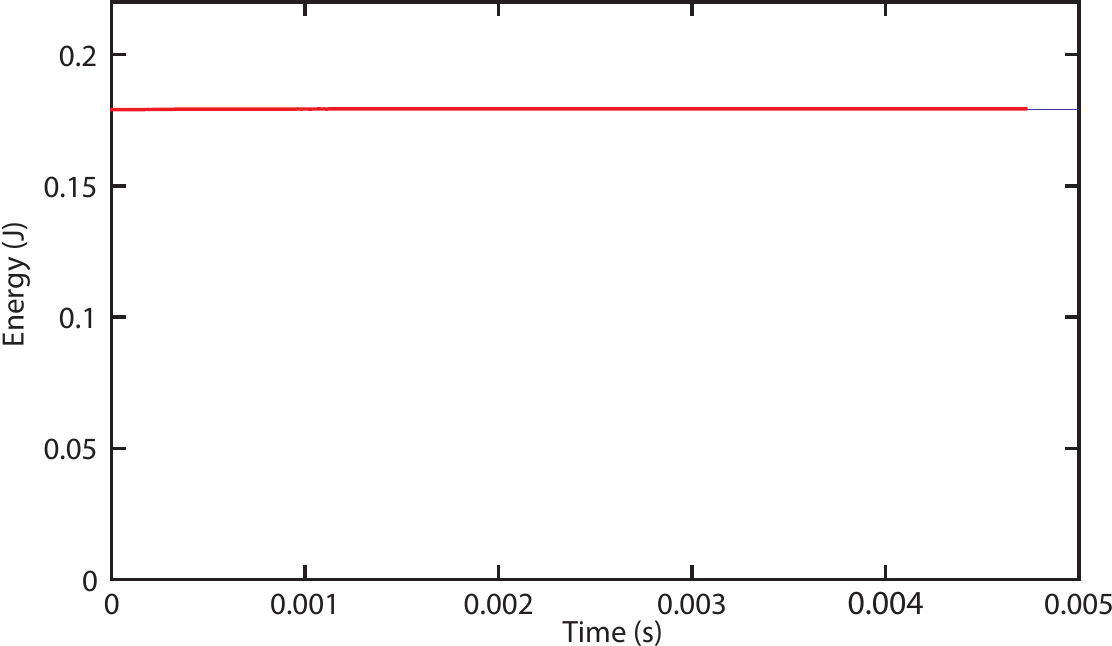}
\caption{Total energy over time of a thin sphere colliding against a
  thin plate, simulated using the proposed contact response method (right) compared to data provided for decomposition
  contact response~\cite{aCirak2005} (left).}
\label{fig-plate}
\end{center}
\end{figure}

Figure \ref{fig-plate} compares energy over time when this simulation is run using both the proposed method and DCR.
Using the former there is no noticeable long-term drift; closely examining the energy data reveals the
high-frequency, low-amplitude, qualitatively-negligible oscillations
characteristic of symplectic integration. The latter introduces noticeable artifical energy decay.

\subsection{Large-scale Three-dimensional Examples}

Harmon \etal~\cite{aHarmon2009} describe a series of optimizations that improve the efficiency of Algorithm \ref{alg-AVIslab}. These optimizations were incorporated to form our Asynchronous Contact Mechanics (ACM) code. This code continues to yield qualitatively good results when scaled to additional large-scale problems. 

\begin{figure}[htb]
\begin{center}
\includegraphics[width=3in]{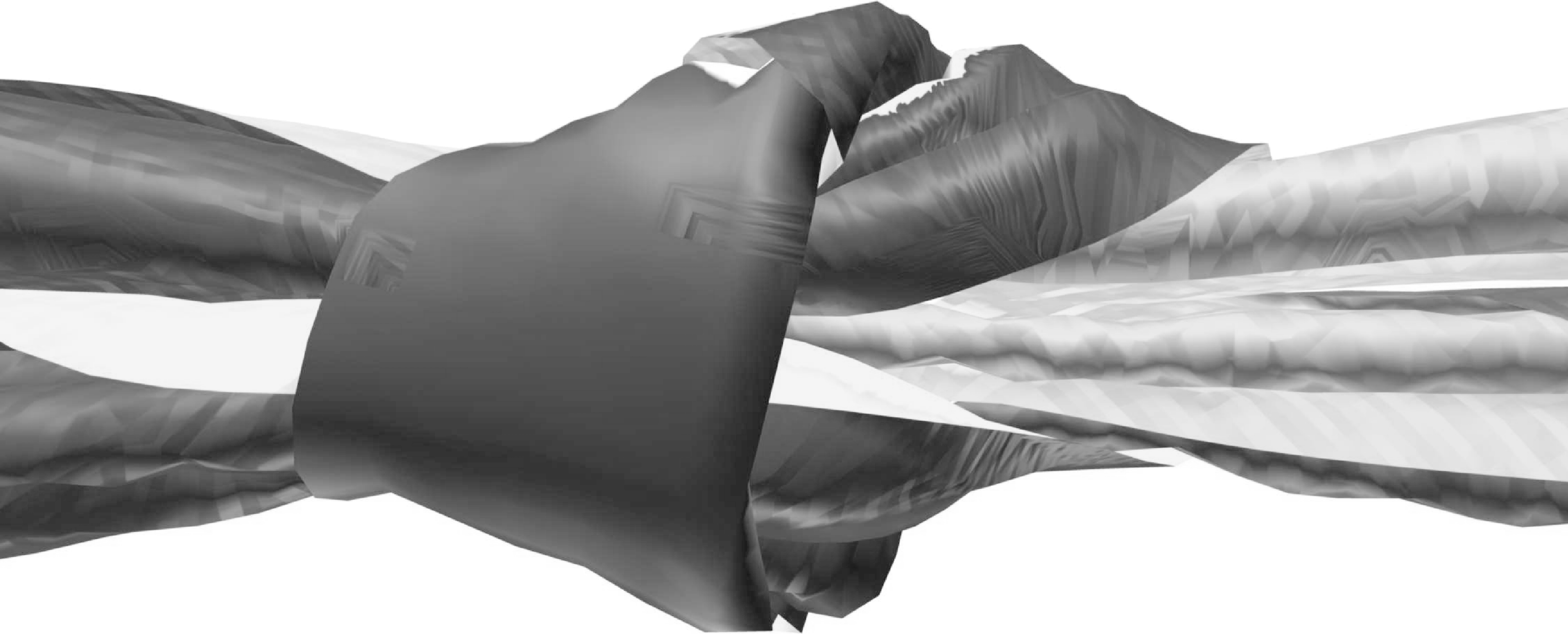}
\includegraphics[width=3in]{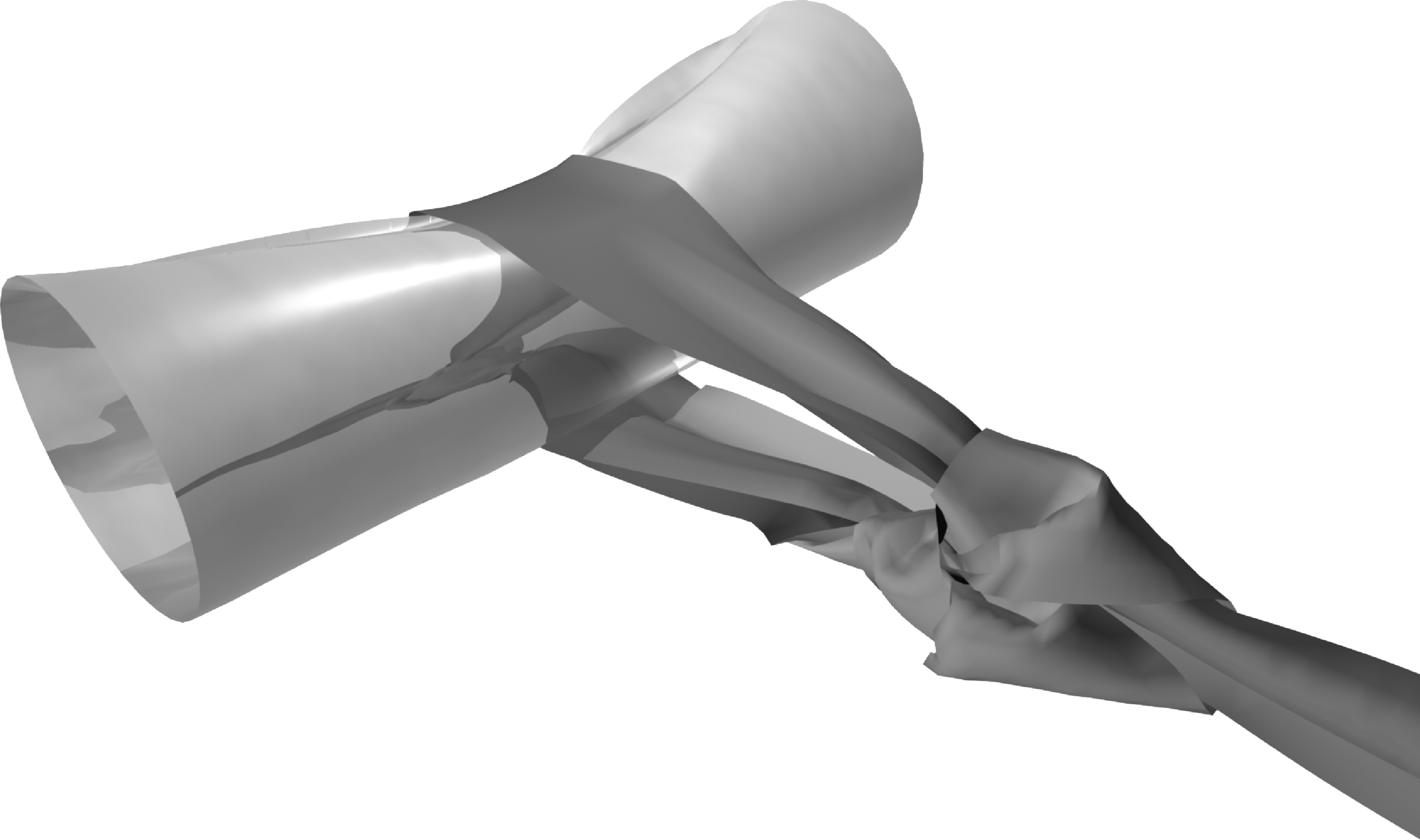}
\end{center}
\caption{Simulated tying of a cloth reef knot (left) and bowline knot (right).}
\label{fig-ribbons}
\end{figure}

Two thin rectangular $27 \unit{cm} \times 2 \unit{cm}$ ribbons were modeled as thin shells of $5321$ vertices, subject to constant-strain triangle stretching forces~\cite{rYotam2004} (stiffness $750$) and discrete shell bending forces formulated by Grinspun \etal~\cite{aGrinspun2003} (stiffness $0.05$). These ribbons were positioned into a loose reef knot by an artist. The knot was then tightened by constraining the end of the ribbon to move apart at $10 \unit{cm/s}$, and running the simulation.

Figure \ref{fig-ribbons}, left, shows the ribbon after $2$ seconds. Since the velocities of the ends of the ribbons were constrained, the knot material became arbitrarily stretched once the knot was tight. The forces pressing the two ribbons into each other thus grew unbounded, but the two ribbons never interpenetrated, nor were other collision-related artifacts observed. It should be stressed that this good behavior did not require the tweaking of the penalty stiffnesses nor any other artificial parameters.

As a second large-scale example, a ribbon similar to the ones in the reef knot simulation was positioned by an artist into a loose bowline knot tied around a cylindrical thin shell of $1334$ vertices. The bowline was then tightened by fixing one end of the ribbon and constraining the other to move away from the cylinder at $10 \unit{cm/s}$. Again, the knot successfully tightened with no penetrations or other artifacts (figure \ref{fig-ribbons}, right).

\subsection{Sphere and Wedge}
\begin{figure}[htb]
\begin{center}
\includegraphics[width=2in]{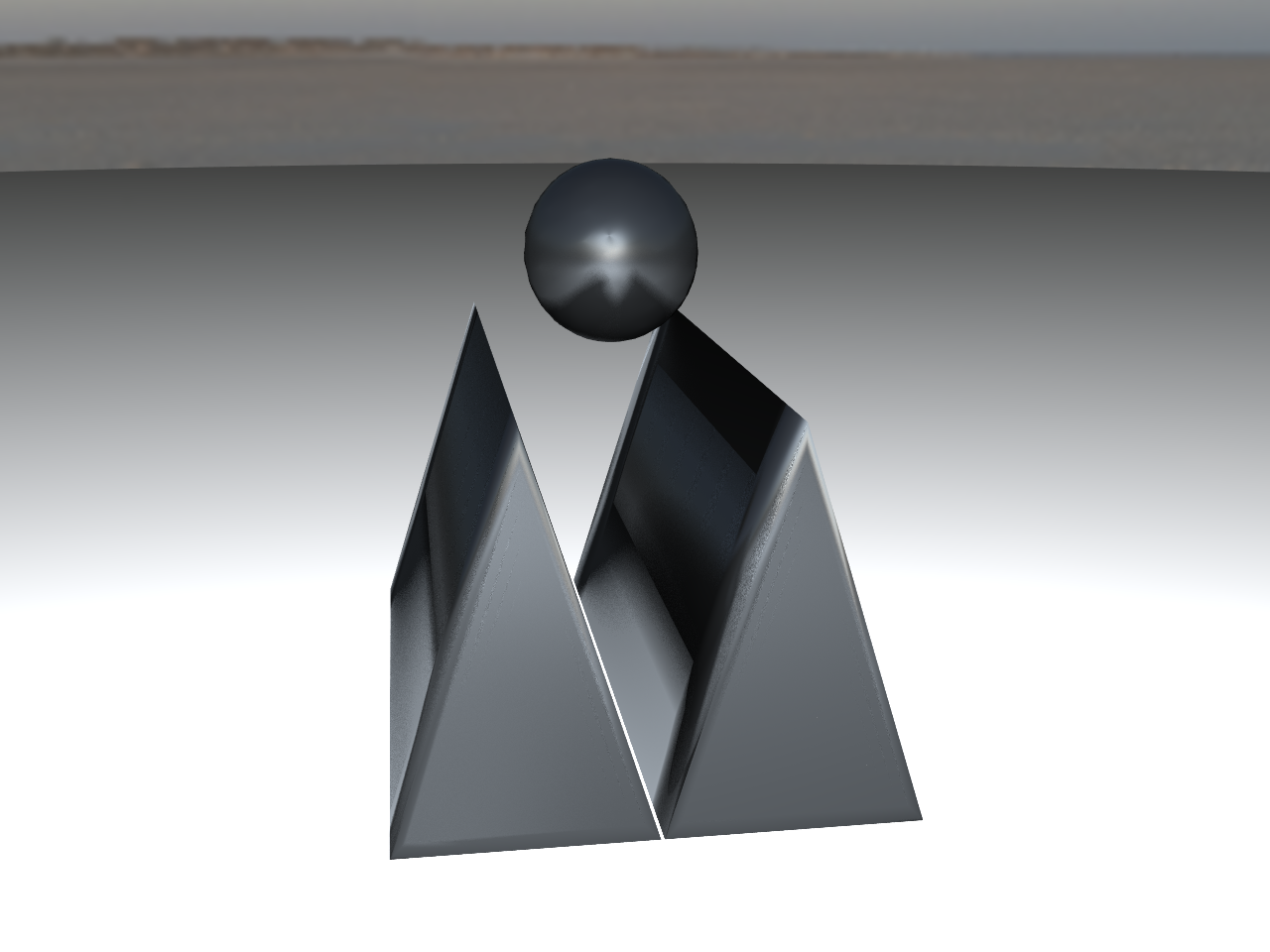}
\includegraphics[width=2in]{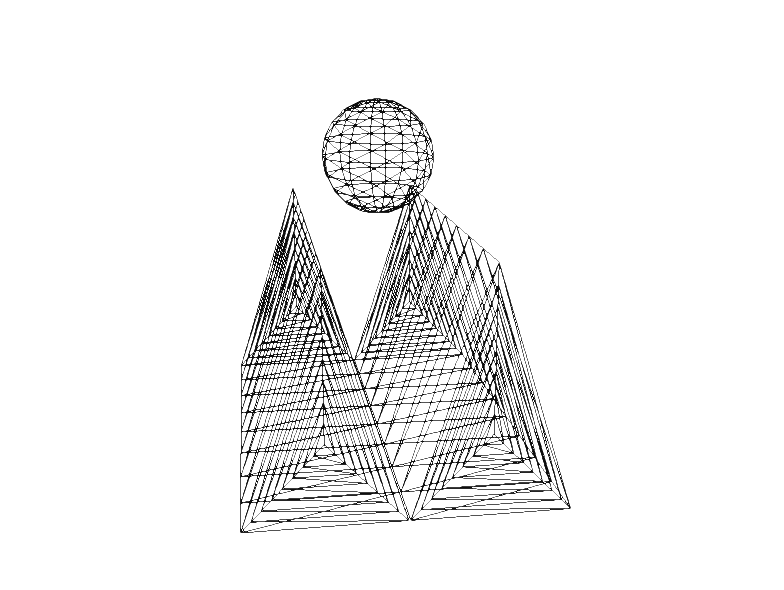}
\includegraphics[width=2in]{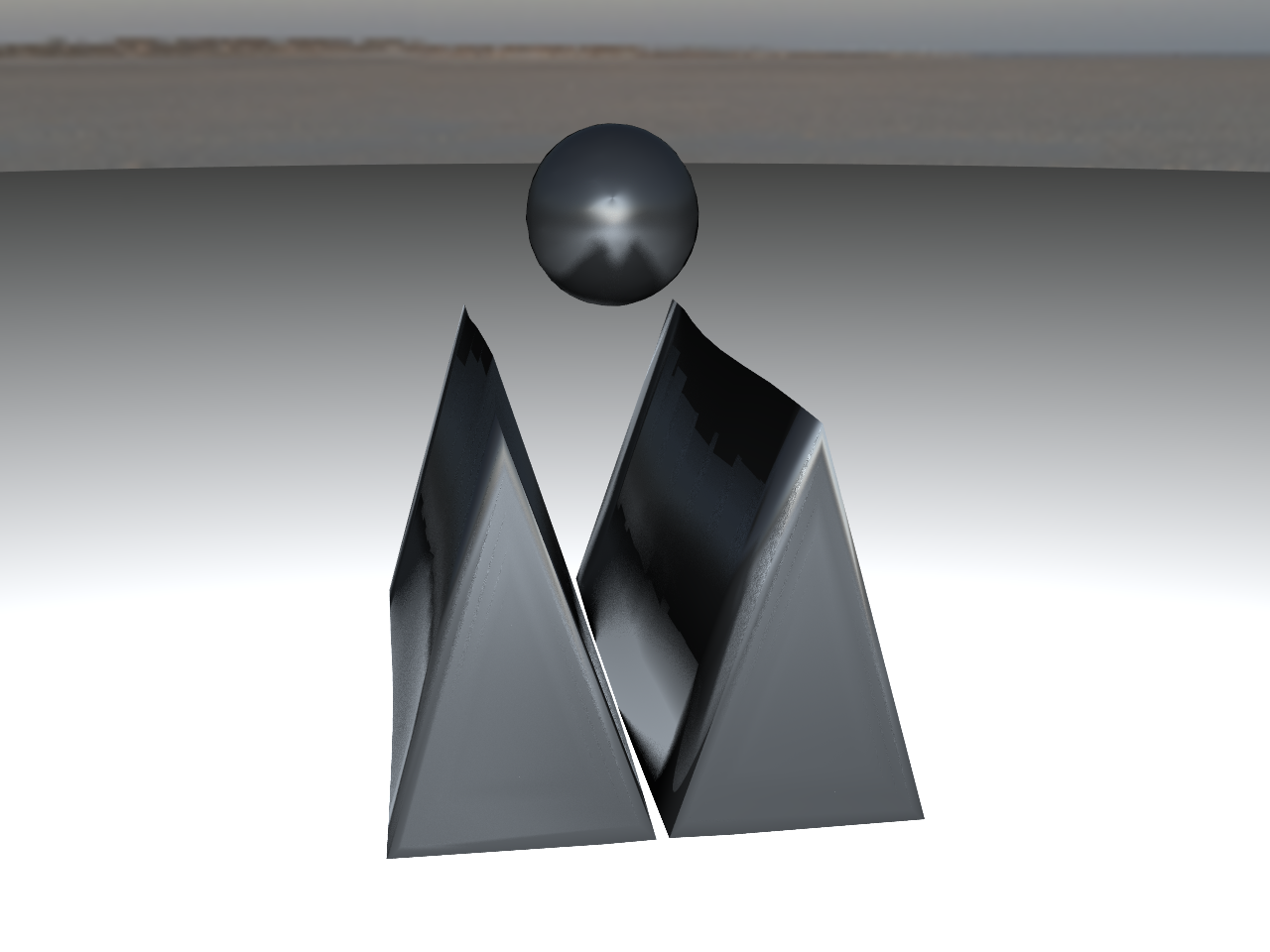}
\end{center}
\caption{A sphere falling into a wedge, at the beginning of the simulation (left and center) and 0.42 seconds later, after the sphere has reflected off of the wedge (right). The center figure shows the mesh elements of the bodies.}
\label{fig-wedge}
\end{figure}

Inspired by Pandolfi \etal~\cite{aPandolfi2002}, a rigid thin-shell sphere was dropped into a wedged formed by two thin shell triangular prism, shown in Figure \ref{fig-wedge}. Each prism has an isosceles base with width $12.92 \unit{cm}$ and height $20.05 \unit{cm}$, and length $38.41 \unit{cm}$. The prisms contain 71 vertices each. The sphere contains 92 vertices, has radius $4.97\unit{cm}$ and begins the simulation $20.84\unit{cm}$ above the ground plane on which the prisms rest. The sphere has initial downwards velocity of $-100\unit{cm/s}$ (no gravity). The sphere and shells use the same thin shell model as the debris in the above trash compactor example, with bending and stretching stiffness parameters $100000$ and $50000$ respectively.

\begin{figure}[!htb]
\begin{center}
\includegraphics[width=3in]{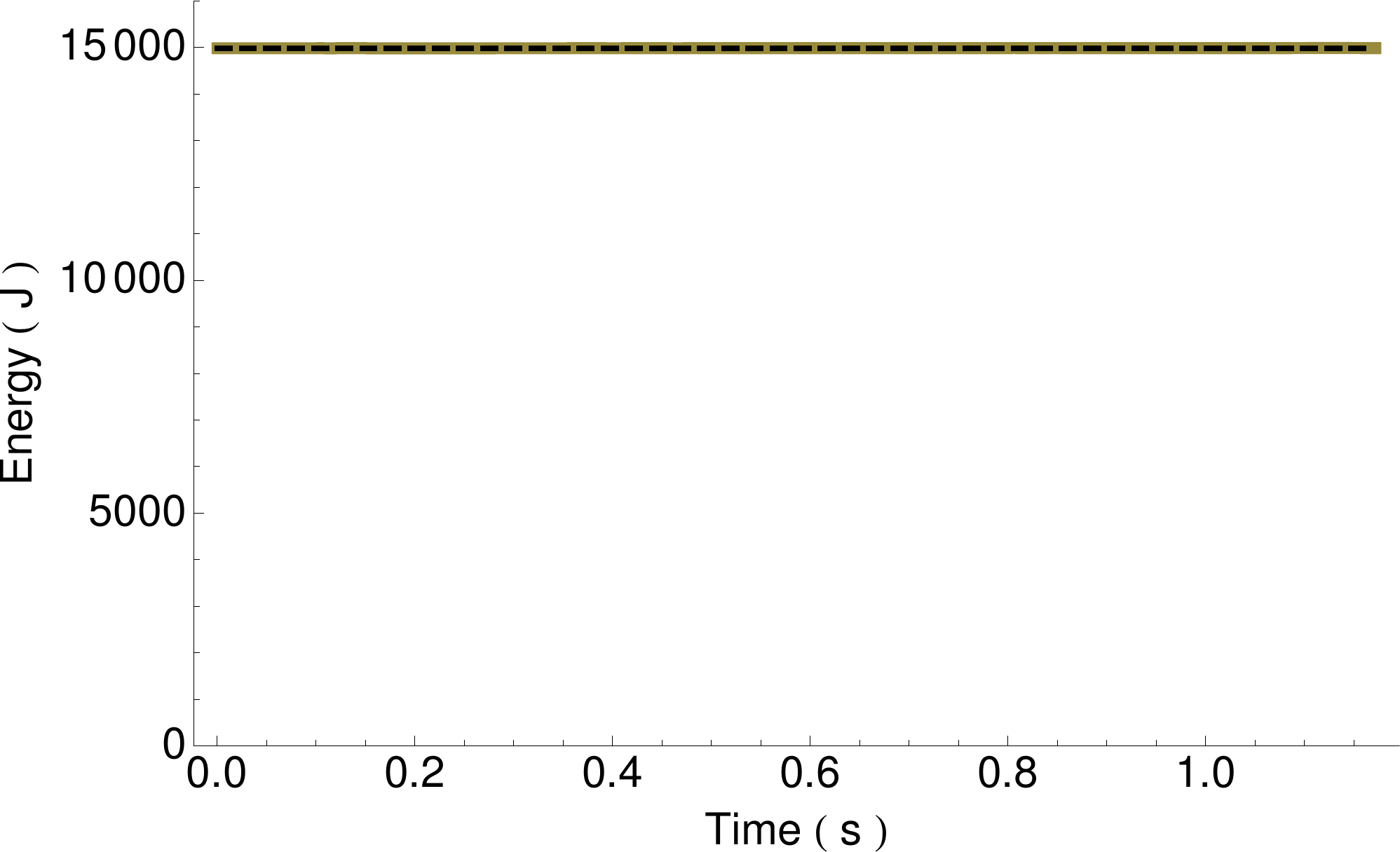}
\end{center}
\caption{The energy of the wedge-sphere system as a function of time. The energy does not drift over the course of the simulation.}
\label{fig-wedgeenergy}
\end{figure}
As the sphere descends, it enters into multiple contact with the faces of the wedge, which undergo elastic deformation and high-frequency vibration. Despite the large areas of simultaneous contact and high velocity at the time of impact, the energy of this system, plotted in Figure \ref{fig-wedgeenergy}, exhibits good behavior and does not drift.

\subsection{Draping on Spikes}

\begin{figure}[!htb]
\begin{center}
\includegraphics[width=3in]{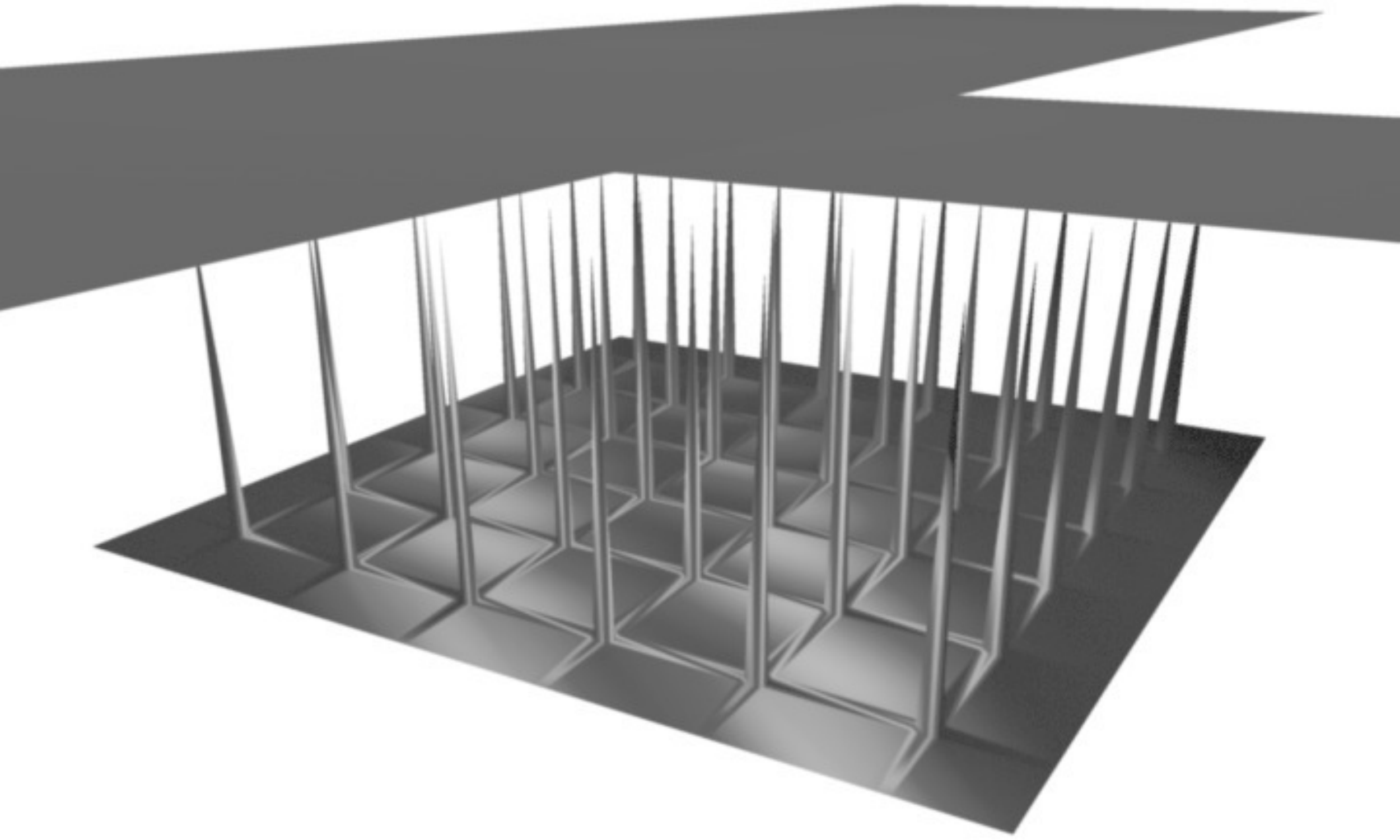}
\includegraphics[width=3in]{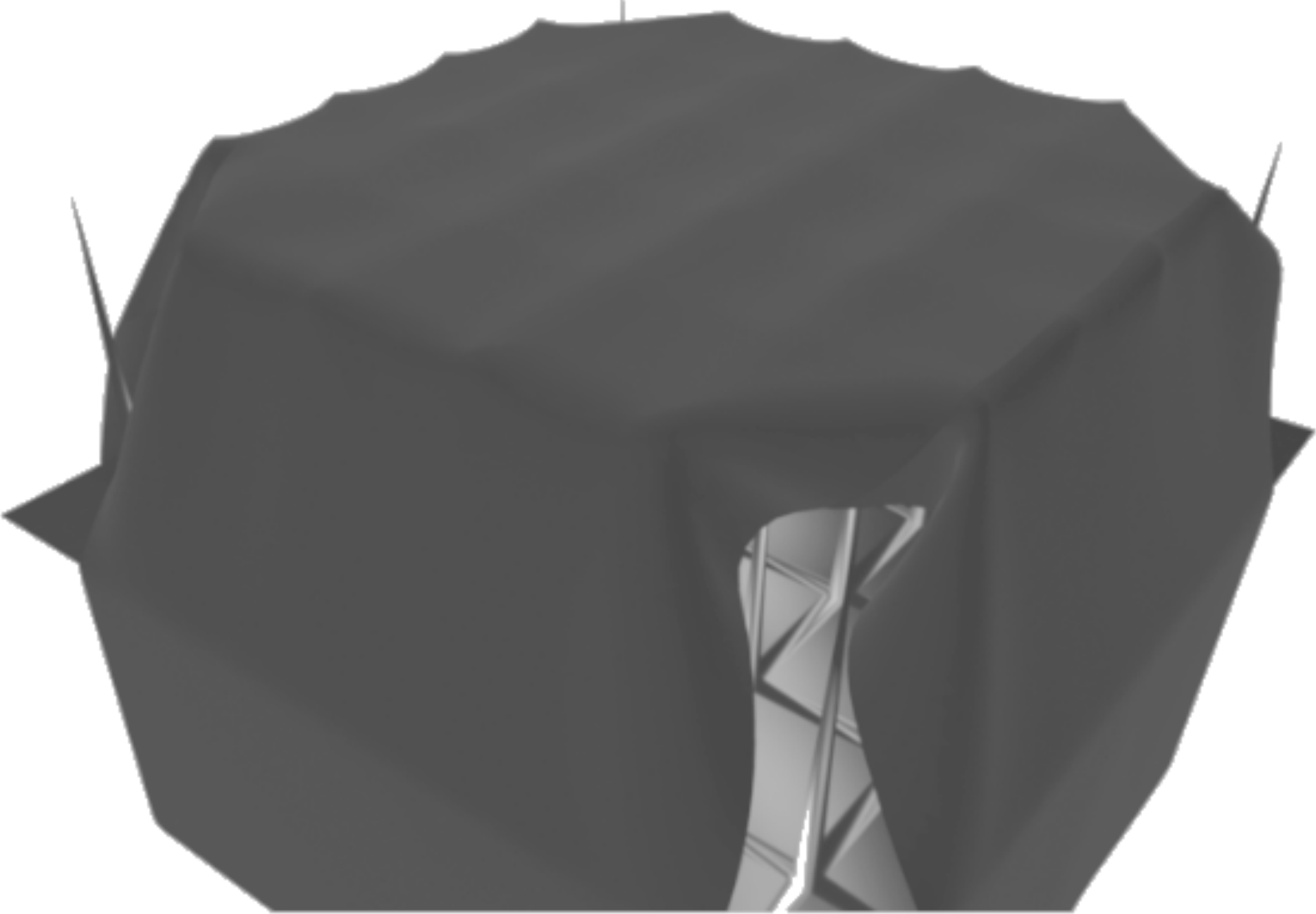}
\end{center}
\caption{Two cloth rectangles were draped on a bed of spikes. The system at the start of the simulation (left), and after the cloth has come to rest (right).}
\label{fig-spikes}
\end{figure}

\begin{figure}[!htb]
\begin{center}
\includegraphics[width=3in]{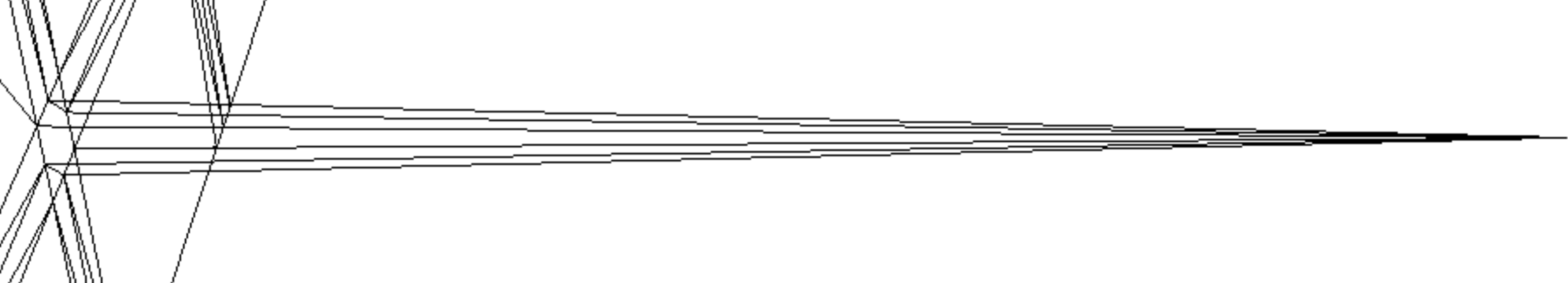}
\end{center}
\caption{A close-up of one of the spikes; the spike has been rotated clockwise 90 degrees to conserve space. Each spike is composed of six triangles with apex angle $3.47$ degrees.}
\label{fig-spikesclose}
\end{figure}

ACM's ability to robustly handle degenerate geometry was tested by dropping two 1994-vertex, $15 \unit{cm} \times 50 \unit{cm}$ cloth meshes (stretching stiffnes $500$, bending stiffness $0.1$, stretching damping $1.0$, bending damping $0.1$) on top of a rigid $20 \unit{cm} \times 20 \unit{cm}$ plate from which protrude 36 $7.8 \unit{cm}$ spikes (see figure \ref{fig-spikes}). Each spike was modeled using six highly-degenerate, sliver triangles: each triangle's most acute angle measures 3.47 degrees (see figure \ref{fig-spikesclose}). The cloth was allowed to fall under gravity ($9.8 \unit{m/s^2}$) and drape on top of the spikes until it had come to rest. No penetrations, oscillations, or other artifacts were observed. 

\begin{figure}[!htb]
\begin{center}
\includegraphics[width=3in]{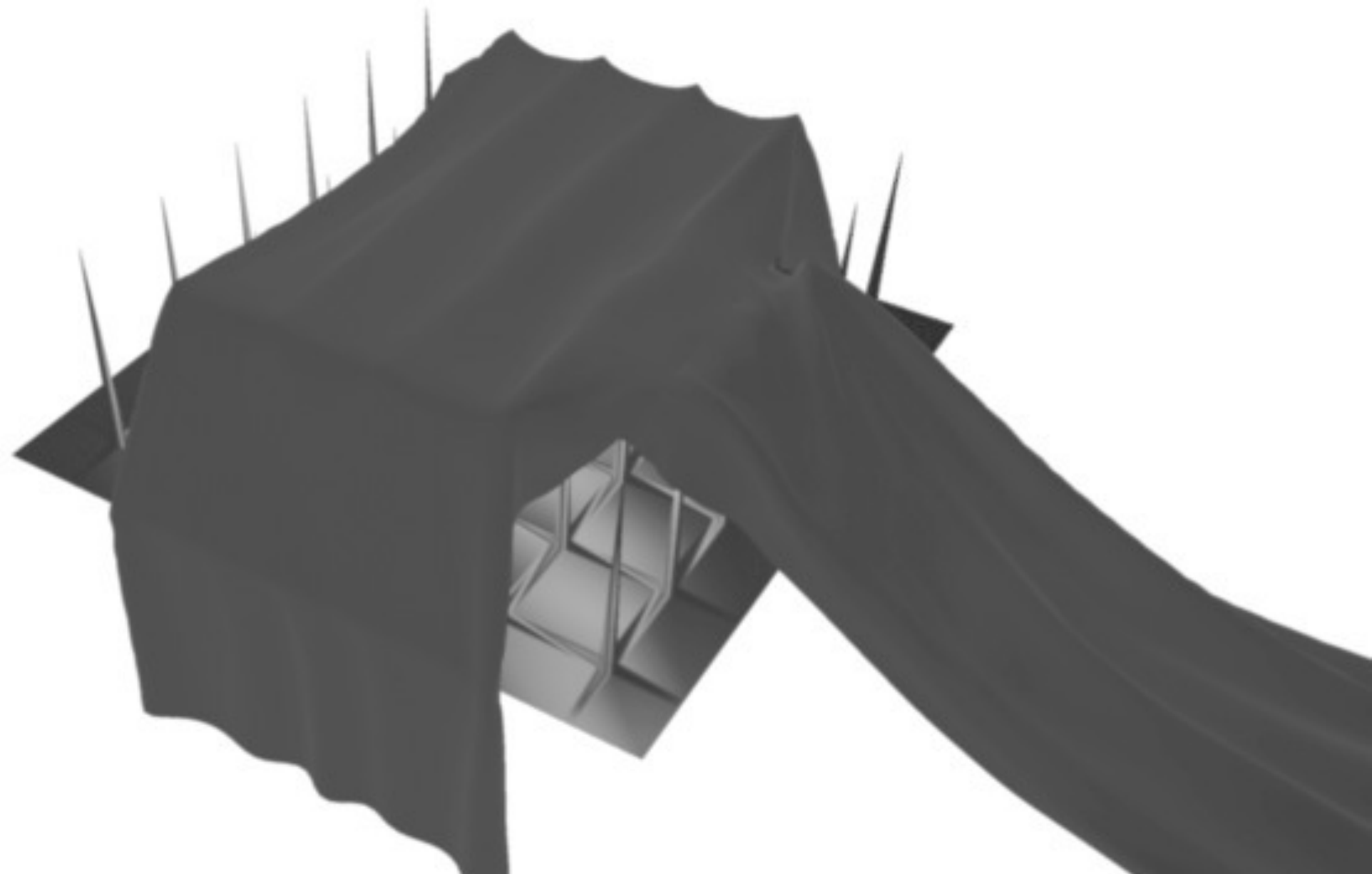}
\end{center}
\caption{After the cloth came to rest, the bottom cloth was pulled out from under the top one. The simulation 3 seconds after pulling began.}
\label{fig-spikespull}
\end{figure}

After the cloth came to rest, the bottom cloth was pulled out from under the top one by constraining one side of the cloth to move at $10 \unit{cm/s}$ parallel to and away from spiked plate; see figure \ref{fig-spikespull}. The bottom cloth scraped against the spikes and slid, with no dissipation, against the top cloth. No interpenetrations occured.

\subsection{Trash Compactor}
\begin{figure}[!htb]
\begin{center}
\includegraphics[width=6in]{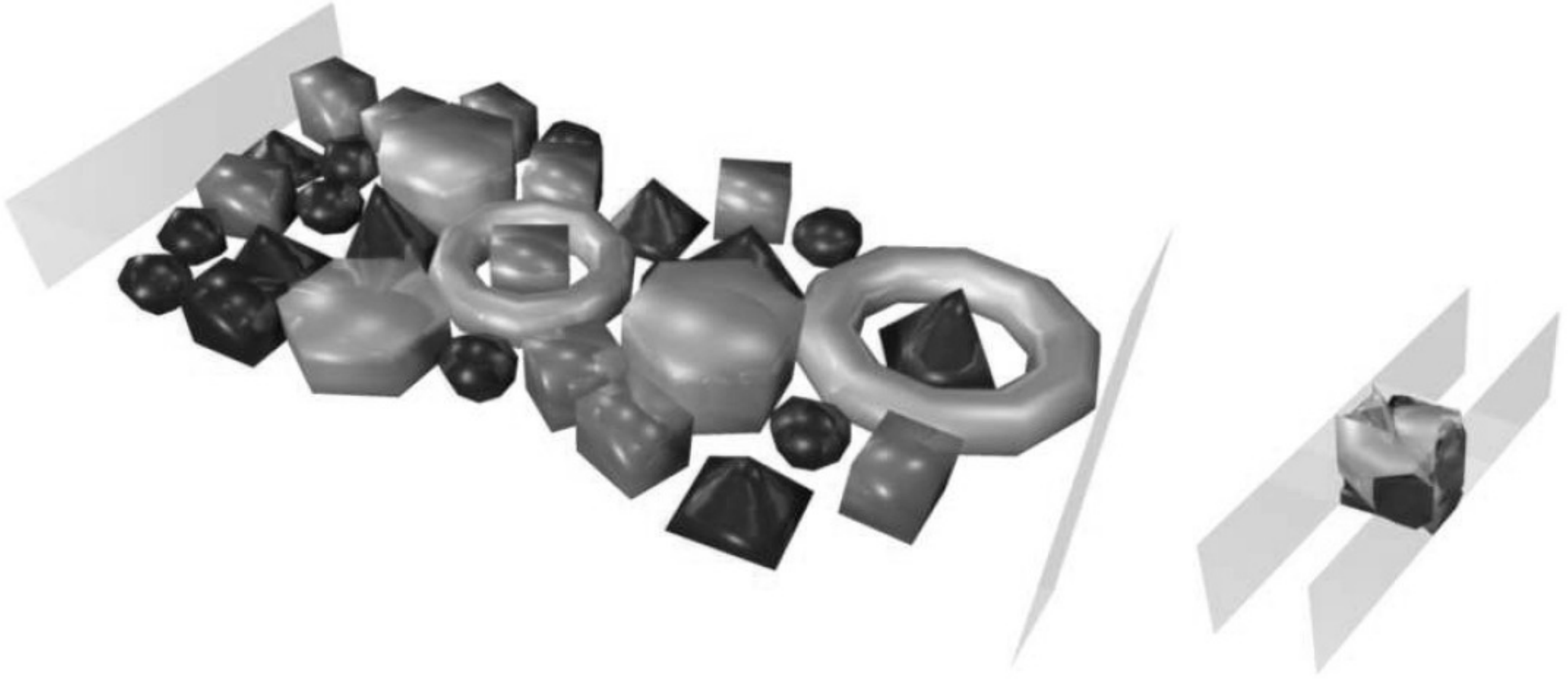}
\end{center}
\caption{Walls close in and compress various thin-shell objects. The beginning (left) and end (right) of the simulation.}
\label{fig-compact}
\end{figure}

Various coarse thin-shell solid objects (platonic solids, tori, etc.) modeled as triangle meshes were placed in a rectangular box measuring $71.5 \unit{cm} \times 36.7 \unit{cm} \times 9.3 \unit{cm}$. The four sides were scripted to close in and compress the objects within: the length at $20 \unit{cm/s}$, and the width at $10 \unit{cm/s}$. All objects were given the same material parameters (stretching stiffness $1000$, bending stiffness $10$, stretching damping $15$, bending damping $0.5$) and held to the floor of the box by gravity ($9.8 \unit{m/s^2}$).  Figure \ref{fig-compact} shows the box at the beginning of the simulation, and after the simulation had run for $3.4$ seconds. A simple plastic deformation model, described by Bergou \etal~\cite{aBergou2007}, allowed the objects to crush plastically when stressed by the encroaching walls. Nevertheless, the material forces acting on the objects grew larger as the box decreases to a small fraction of its original volume, yet no object penetrated any other object or wall, as guaranteed by the method.

\section{Conclusion and Future Work}
A framework for asynchronous, structure-preserving handling of contact and impact has been presented. Provable guarantees were established for this framework: impact handling is robust, allowing no penetrations or tunneling; the good properties of AVIs are preserved, such as a discrete Noether's Theorem and discrete multisymplectic structure; and for well-posed problems, the amount of computation required to simulate the problem is bounded and in particular finite. Good long-time energy behavior, conjectured to accompany multisymplectic structure, was confirmed empirically by both didactic and challenging, large-scale experiments. Modifications to allow for simple dissipative phenomena, such as a coefficient of restitution, were described. Data structures and algorithms to improve performance, such as the use of separating slabs to prune inactive penalty layers, were briefly discussed. Implementation details for these and other optimizations, as well as ideas for future improvements to the algorithms that promise to substantially decrease computation time, are treated more comprehensively by Harmon \etal~\cite{aHarmon2009}. 

Missing from the basic asynchronous contact framework described by this paper is comprehensive handling of friction, particularly static friction. Static friction conflicts fundamentally with asynchrony: in an asynchronous simulation, contact between a pair of elements is resolved piecemeal, by summing the impulses at many different times contributed by many different penalty layers. At any given moment of time it is unclear how to define a total normal force, an element necessary for the robust treatment of even the most elementary static friction models. Successfully merging the handling of friction with the asynchronous framework, to allow simulations of systems such as a standing house of cards, remains a challenging area for future research.

\section{Acknowledgements}
The authors thank David Mooy for modeling the
knots, and Igor Boshoer, Matt Kushner, Kori Valz for lighting and
rendering. We are grateful for the valuable feedback provided by
Breannan Smith, Mikl\'{o}s Bergou, Rony Goldenthal, Bernhard Thomaszewski, Max
Wardetzky, and the anonymous reviewers. This work was supported
in part by the Sloan Foundation,the NSF (MSPA Award Nos. IIS-05-28402 and IIS-09-16129, CSR Award
No. CNS-06-14770, CAREER Award No. CCF-06-43268), and the
Amazon Elastic Compute Cloud. The Columbia authors are supported
in part by generous gifts from Adobe, ATI, Autodesk, mental
images, NVIDIA, Side Effects Software, the Walt Disney Company, and Weta Digital.

\clearpage
\bibliography{refs}
\bibliographystyle{model1-num-names}
\end{document}